

\input amstex
\documentstyle{amsppt}

\magnification=\magstep1
\TagsOnLeft
\NoBlackBoxes 
\pageheight{9truein}
\pagewidth{6.3truein}

\def \Barnes{1}
\def \Barnestwo{2}
\def \Cardon{3}
\def \Lemm{4}

\def \CusickFlahive{6}
\def \Godwin{5}
\def \paperone{6}
\def \Var{7}
\def\sym{{\sim}}
\def\ve{{\varepsilon}}

\define\Z {\bold Z}
\define\R {\bold R}
\define\N {\bold N}
\define\Zp {\bold Z_p}
\define\Q {\bold Q}


\topmatter
\title  On the inhomogeneous spectrum of
period two quadratics \endtitle 
\leftheadtext{CHRISTOPHER G. PINNER}
\author{ Christopher G. Pinner  }\endauthor
\rightheadtext{ On the inhomogeneous spectrum  }
 
\abstract{For an  $\alpha$ with period two negative
continued fraction expansion:
$$ \alpha =[0;a_1,a_2,a_3,...,a_N,\overline{a,b}]^-=\frac{1}{\displaystyle a_1-\frac{1}{\displaystyle a_2-\frac{1}{\displaystyle a_3-\cdots}}},\;\;\;2\leq a<b, $$
we give a complete description of the values in its inhomogeneous
Lagrange spectrum
$$ {\bold L}(\alpha):=\{M(\alpha,\gamma) \; : \; \gamma \not\in \Z+\alpha \Z\}, \;\; M(\alpha,\gamma):=\liminf_{n\rightarrow\infty} |n| ||n\alpha-\gamma|| $$
down to the first limit point. The largest point in the spectrum is
always isolated (as conjectured by Barnes and Swinnerton-Dyer).

If $a$ is odd then   
the first limit point is usually the second or third largest value (the fourth in two cases $(a,b)=(5,7)$ or $(7,9)$, and a limit from above when $(a,b)=(3,4),(3,5),(3,6)$). When $a$ is even and $b$ is odd there are infinitely many values above the first limit point.
This is also the case when $a$ and $b$ are both even except
for $b=2a-2$, $b=a+2$ or $b=a+4$, $a\geq 4$,  when the limit point is usually the third largest value. } \endabstract
 
 
 
\address{Mathematics and Computer Science,
University of Northern British Columbia,
3333 University Way, Prince George, B.C., Canada V2N 4Z9.
{\tt pinner\@math.ksu.edu}}
\endaddress

\endtopmatter
\document
\baselineskip=18pt 

\head 1. Introduction. \endhead

For a pair of real numbers $\alpha$ and $\gamma$, with $\alpha$ irrational,
one defines the inhomogeneous approximation constant
$$  M(\alpha,\gamma):=\liminf_{n\rightarrow\infty} |n| ||n\alpha-\gamma||. $$
Of course $M(\alpha,\gamma)=M(\alpha,n\alpha+m\pm \gamma)$ for any $n,m$ in $\Z$. Fixing $\gamma=0$ and varying $\alpha$ gives the (reciprocal of the) traditional Lagrange spectrum 
$${\bold L}^{-1} =\{M(\alpha,0)\; :\; \alpha \in \R\}.  $$
 Here we instead fix $\alpha$ and vary $\gamma$ 
to obtain  the inhomogeneous Lagrange
spectrum for $\alpha$
$$   {\bold L}(\alpha):=\{M(\alpha,\gamma) \; : \; \gamma \not\in \Z+\alpha \Z\}. $$
For the classical Lagrange spectrum old results of Markoff and Lagrange
(see [\CusickFlahive]) give a complete description of  ${\bold L}^{-1}$
from the largest value $1/\sqrt{5}$ down  to the first limit point $1/3$.
We suppose here  that 
$$ \alpha=[0;a_1,...,a_{N-1},\overline{a,b}]^-=\frac{1}{\displaystyle a_1-\frac{1}{\displaystyle a_2-\frac{1}{\displaystyle a_3-\cdots}}}, \;\;\; 2\leq a<b, \tag1.1 $$
has a period two negative expansion and give a similarly complete description of
${\bold L}(\alpha)$ above the first limit point. We dealt with the period one 
cases $a=b$  in [\paperone]. The first limit point is always a limit from below but need not be a limit point from above, indeed it can even be the second largest
point, as was observed for the related inhomogeneous Markoff 
spectrum [\Godwin].
Godwin's example $x^2-23y^2$ of a form with non-isolated second minimum 
corresponds to the special case  $a=5$, $b=10$ ($r=10,\; m=0$ in Theorem 2).
It has been conjectured by Barnes--Swinnerton-Dyer [\Barnes] that the 
largest value in ${\bold L}(\alpha)$, denoted $\rho (\alpha)$,  is isolated for any quadratic $\alpha$.
We show that this is certainly true for all period two $\alpha$.
Many special cases of $\rho (\alpha)$ occur in [\Barnes], including the
the period two examples $(a,b)=(3,3n)$ and $(a,b)=(2n+1,2(2n+1))$.
We recall (see [\Barnes]) that  $\Z[\alpha]$ is norm Euclidean 
with respect to  $(x-\alpha y)(x-\bar{\alpha} y)=x^2+bxy+cy^2$ 
iff $\rho (\alpha)< 1/\sqrt{b^2-4c}$  (see also [\Barnestwo,\Lemm]).
It is interesting to note that $\sqrt{14}=[4;\overline{4,8}]^-$
(where $\Q(\sqrt{14})$ is a field recently shown by Malcolm Harper  to be be Euclidean 
but not norm-Euclidean (see also [\Cardon]))
has just one point in its spectrum above $\frac{1}{2\sqrt{14}}$.

Setting
$$ \align \eta   & :=[0;\overline{a,b}]^-=\frac{1}{2}(b-\sqrt{b^2-4b/a}), \\
\;\;\; \beta   & :=[0;\overline{b,a}]^-=\frac{1}{2}(a-\sqrt{a^2-4a/b}),
\endalign $$
and
$$  D:=\eta \beta, $$
we  define 
$$    \;\;\; M^*(\alpha,\gamma):=4(1-D)M(\alpha,\gamma). $$ 
We write $\rho^*(\alpha)=\rho_1^* (\alpha),\rho_2^* (\alpha), ...$ for
the successive values  of $M^*(\alpha,\gamma)$ above the first 
accumulation point.
For a given $\alpha$ we showed in [\paperone] how to use an $\alpha$-expansion
of $\gamma$ to evaluate $M(\alpha,\gamma)$. Since ${\bold L}(\alpha)={\bold L}(\eta)={\bold L}(\beta)$ it will be convenient to assume throughout  that 
$\alpha$ is purely periodic with 
$a_{2i-1}=a$ and $a_{2i}=b$ for $i\geq 1$, 
so that the $\alpha$-expansions take the form
$$ \gamma =\sum_{i=1}^{\infty} (b_{2i-1}\eta +b_{2i}D)D^{i-1}. $$
Defining integers $t_i$ by $b_i=\frac{1}{2}(a_i-2+t_i)$ and setting
$$\align d_i^+ & :=\sum_{j=0}^{\infty} \left(t_{i+2j+1}\alpha_i +t_{i+2j+2}D\right) D^{j}, \\
d_i^- & :=\sum_{j=0}^{\infty} \left(t_{i-2j}\alpha_{i-1} +t_{i-2j-1}D\right) D^{j}, \endalign $$
where 
$$ \alpha_j:=\cases \beta & \hbox{ if $j$ is odd,}  \\
\eta & \hbox{ if $j$ is even,}  \endcases $$
we showed that if $t_k=a_k$ at most finitely often then
$$ M^*(\alpha,\gamma)=
\liminf_{k\rightarrow \infty} \min\{ s_1^*(k),s_2^*(k),s_3^*(k),s_4^*(k)\} $$
with
$$ \align 
s_1^*(i)  & = (1-\alpha_i+d_i^+)(1-\alpha_{i-1}+d_{i}^-), \\
s_2^*(i)  & = (1+\alpha_i-d_i^+)(1+\alpha_{i-1}+d_{i}^-), \\
s_3^*(i)  & = (1-\alpha_i-d_i^+)(1-\alpha_{i-1}-d_{i}^-), \\
s_4^*(i)  & = (1+\alpha_i+d_i^+)(1+\alpha_{i-1}-d_{i}^-). \endalign $$
When $t_k=a_k$ infinitely often then we need to check the minimum
of $s_1^*(k)$ and $s_2^*(k)$ for both both $\gamma$ and $1-\gamma$. When 
$t_k\neq a_k$ the expansion of $1-\alpha-\gamma$ can be  obtained from
that of  $\gamma$ by simply replacing the
$t_i$ with $-t_i$. Changing a finite number of $t_i$ of course does not affect $M(\alpha,\gamma)$.
We also showed  in [\paperone] 
that if $t_k=a_k$ infinitely often with $k\equiv j$ mod 2
then 
$$ M^*(\alpha,\gamma) \leq \alpha_{j-1}   \tag 1.2 $$
(so we only need worry about these when $a_j$ is small)
and if $t_i= a_i$ at most finitely often with $|t_k|\geq t$, 
$k\equiv j$ mod 2 infinitely often then
$$  M^*(\alpha,\gamma) \leq (a_j-t) \alpha_{j-1}=1-t\alpha_{j-1}+D.\tag 1.3 $$
In view of this last result we would expect the $t_i$ 
corresponding to the larger values of $M(\alpha,\gamma)$ to be all small 
(although when the smallest partial quotient $a$ is odd we in fact
typically  compensate
for the impossibilty of $t_{2i-1}=0$ by selecting the $t_{2i}$ to be roughly $b/a$ away from zero).

We shall define sets $S_k$ of $\gamma$ with eventually periodic $\alpha$-expansions $b_i$. In particular the $\gamma$ of interest will have their
expansions built from the blocks $(b_{2k-1},b_{2k})$ of the form
$$\align
 A_t= \left(\frac{1}{2}(a-2),\frac{1}{2}(b-2+t)\right),\hskip8ex &  {A'}_{t}= \left(\frac{1}{2}(a-2),\frac{1}{2}(b-2-t)\right), \\
 B_t= \left(\frac{1}{2}(a-2-1),\frac{1}{2}(b-2+t)\right),\hskip3ex & {B'}_{t}=\left(\frac{1}{2}(a-2+1),\frac{1}{2}(b-2-t)\right), \\
 C_t= \left(\frac{1}{2}(a-2-2),\frac{1}{2}(b-2+t)\right),\hskip3ex &  {C'}_{t}= \left(\frac{1}{2}(a-2+2),\frac{1}{2}(b-2-t)\right),\\ 
 E_t= \left(\frac{1}{2}(a-2-3),\frac{1}{2}(b-2+t)\right),\hskip3ex & {E'}_{t}=\left(\frac{1}{2}(a-2+3),\frac{1}{2}(b-2-t)\right), \\
 F_t= \left((a-1),\frac{1}{2}(b-2-t)\right),\hskip3ex & {F'}_{t}=\left(
(a-1),\frac{1}{2}(b-2+(t-4))\right). \\
\endalign $$

\head 2. A complete description  down to the first limit point.  \endhead

\head The spectrum when $a\geq 4$ is even. \endhead

Suppose first that $a\geq 4$ is even.

When $b$ is odd we define the classes of $\gamma$
$$\align
S_{-1} & =\{\gamma \; :\; \hbox{ $\gamma$ has periodic expansion, period $A_1A_1A_1'A_1'$}\}, \\
S_{-2} & =\{\gamma \; :\; \hbox{ $\gamma$ has periodic expansion, period $C_3$}\}, 
\endalign $$
and for integers $k\geq 0$
$$S_k  =\{\gamma \; :\; \hbox{ $\gamma$ has periodic expansion, period $A_1'A_1 (A_1'A_1'A_1A_1)^k$}\}. $$
For $b$ even we set
$$ \align
S_{k,1} & =\{\gamma : \hbox{ $\gamma$ has periodic expansion, period $A_0 (A_2{A'}_{2})^k$} \},\;\; k\geq 0,\\
S_{k,2} & =\{\gamma : \hbox{ $\gamma$ has periodic expansion, 
period $A_2C_4^kC_2{A'}_2{C'}_{4}^k{C'}_2$} \}, \;\;k\geq 1, \\
S_{k,3} & =\{\gamma : \hbox{ $\gamma$ has periodic expansion, period $C_4 (C_2C_4)^k$} \},\;\; k\geq 1, \\
S_{k,4} & =\{\gamma : \hbox{ $\gamma$ has periodic expansion, period $C_4 C_2^k$} \},\;\; k\geq 0, \\
S_{k,5} & =\{\gamma : \hbox{ $\gamma$ has periodic expansion, period $A_2C_2^k{A'}_2{C'}_2^k$}\},\;\; k\geq 0,  \\
S_{-1} & =\{\gamma : \hbox{ $\gamma$ has periodic expansion, period $C_2$}\},\;\;  
\endalign $$
and for $(a,b)=(4,6)$ 
$$ S_{-2} =\{\gamma : \hbox{ $\gamma$ has periodic expansion, period $A_2C_2$}\},\;\; $$
for $(a,b)=(8,12)$ 
$$S_{k,6}  =\{\gamma : \hbox{ $\gamma$ has periodic expansion, period $(A_2'C_2'A_2C_2)^k A_2'C_2'C_2'A_2C_2C_2$}\},\;\; k\geq 0,  $$
and for $(a,b)=(6,10)$
$$S_{k,7}  =\{\gamma : \hbox{ $\gamma$ has periodic expansion, period $(C_4C_2)^k A_2'A_2$}\},\;\; k\geq 1,  $$
We write
$$ \delta_{\square}= M^*(\alpha,\gamma),\;\;\gamma \in S_{\square}. $$

\proclaim{Theorem 1}
Suppose that $\alpha=[0;\overline{a,b}]$ with $a<b$  and  $a\geq 4$ even.

If $b$ is odd then
$$ \rho^* (\alpha)=\cases \delta_0, & \hbox{ if $b=a+1$, $a+3$ or $b\geq 2a-3$,}\\
\delta_{-2}, & \hbox{ if $a+5\leq b\leq 2a-5$,}\endcases $$
and the only values of $M^*(\alpha,\gamma)$ greater than 
$$ \delta_{\infty}:=\left(1-\beta -\frac{\beta (1-D)(1+2D^2)}{1+D^2}\right)\left(1-\eta -\frac{D(1-D)}{1+D^2}\right)$$
are $\delta_{-1}$, $\delta_{-2}$ (when $a+3\leq b\leq 2a-3$) 
and the $\delta_{k}\searrow \delta_{\infty}$, $k\geq 0$.

When $b$ is even
$$ \rho^* (\alpha)=\delta_{0,1}=(1-\eta)(1-\beta ). $$
For $b\geq 2a$ or $(a,b)=(4,6)$ the values of $M^*(\alpha,\gamma)$ greater than
$$ \delta_{\infty,1}=\left(1-\eta +\frac{2D^2}{1+D}\right)\left(1-\beta-\frac{2\beta }{1+D}\right) $$
are the $\delta_{k,1}\searrow \delta_{\infty,1}$, $k\geq 0$ and $\delta_{0,5}$ plus  $\delta_{0,4}$
when $2a\leq b\leq 3a-6$, and $\delta_{-2}$ when $(a,b)=(4,6)$.

When $b=2a-2$, $a\geq 8$ the 
values of $M^*(\alpha,\gamma)$ greater than
$$ \delta_{\infty,2}=\left(1-3\eta +\frac{2D(2-\eta)}{1-D}\right)\left(1+\beta -\frac{2\beta  D(1-\eta+D)}{1-D}\right) $$
are $\delta_{0,1}$ and $\delta_{0,4}$
(with the $\delta_{k,2}\nearrow \delta_{\infty,2}$).

When $b=2a-4$, $a\geq 10$ the values of $M^*(\alpha,\gamma)$ greater than
$$ \delta_{\infty,3}=\left(1-\eta -\frac{2\eta D}{1-D} +\frac{2D(1+2D)}{1-D^2}\right)\left(1-3\beta +\frac{2D}{1-D} -\frac{2\beta  D(2+D)}{1-D^2}\right) $$
are $\delta_{0,1}$, $\delta_{0,4}$, $\delta_{1,4}$, $\delta_{1,5}$ and the $\delta_{k,3}\searrow \delta_{\infty,3}$, $k\geq 1$.

When $b\leq 2a-6$ or $a=6, b=8$, the first limit point is
$$ \delta_{4,\infty} =\left( 1-\eta +\frac{2D(1-\eta)}{1-D}\right)\left(1-3\beta +\frac{2D(1-\beta )}{1-D}\right) $$
with $\delta_{k,5}\nearrow \delta_{\infty,4}$. The values above 
$\delta_{\infty,4}$
are $\delta_{0,1}$, $\delta_{-1}$, and when $a+6\leq b\leq 2a-6$ the $\delta_{k,4}\searrow \delta_{\infty,4}$.

For $(a,b)=(8,12)$ the values above $\delta_{\infty,6}=\frac{112387809}{2209}\sqrt{138}-\frac{1320256308}{2209}$ are $\delta_{0,1}$, $\delta_{1,5}$ and the $\delta_{k,6}\searrow \delta_{\infty,6}$, $k\geq 0$.

For $(a,b)=(6,10)$ the values above $\delta_{\infty,7}=\frac{6329319}{40}-\frac{6551443}{600}\sqrt{210}$
are $\delta_{0,1}$, $\delta_{0,5}$, $\delta_{1,4}$ with the $\delta_{k,7}\nearrow \delta_{\infty,7}$.

\endproclaim

The values of the $\delta_{\square}$ are given in Lemma 1 below and are achieved only
with $\pm \gamma$ in $S_{\square}$. We
should remark that the limit points are limit points of limit points from below.

We give the proof of Theorem 1 in $\S4$ below.

\head  The spectrum when $a$ is odd.   \endhead  

Suppose now that $a$ is odd.
We suppose that
$$ b=ma+r, \;\;\;0<r\leq 2a,\;\;\; 2\mid r, $$
and set
$$ n:=m+2,\;\;\;\;  s:=m-2,$$
so that $m$ and $n$ are the closest integers to $b/a$ with the same parity as $b$. 
 
We set
$$ \align
S_0 & =\{\gamma : \hbox{ $\gamma$ has periodic expansion, period $B_m$} \}, \\
S_{-1} & =\{\gamma : \hbox{ $\gamma$ has periodic expansion, period $B_mB_n$} \}, \\
S_{-2} & =\{\gamma : \hbox{ $\gamma$ has periodic expansion, period $B_n$} \}, \\
S_{-3} & =\{\gamma : \hbox{ $\gamma$ has periodic expansion, period $B_nB_s$} \}, \\
S_{-4} & =\{\gamma : \hbox{ $\gamma$ has periodic expansion, period $B_mB_m'$ (or $B_sB_s'$ if $m=1$)} \}, \\
S_{-5} & =\{\gamma : \hbox{ $\gamma$ has periodic expansion, period $E_3$ } \}, \\
S_{-6} & =\{\gamma : \hbox{ $\gamma$ has periodic expansion, period $F_2$ } \}, \\
S_{-7} & =\{\gamma : \hbox{ $\gamma$ has periodic expansion, period $F_1$ or $F_3$ } \}, \\
S_{k,1} & =\{\gamma : \hbox{ $\gamma$ has periodic expansion, period $B_n^k B_m$} \}, \\
S_{k,2} & =\{\gamma : \hbox{ $\gamma$ has periodic expansion, period $B_n^k B_m {(B_n')}^{k}B_m'$} \}, \\
S_{k,3} & =\{\gamma : \hbox{ $\gamma$ has periodic expansion, period $(B_mB_n)^kB_n$} \}, \\
S_{k,4} & =\{\gamma : \hbox{ $\gamma$ has periodic expansion, period $(B_nB_m)^k(B_n'B_m')^k$} \}, \\
S_{k,5} & =\{\gamma : \hbox{ $\gamma$ has periodic expansion, period $B_m^kB_n$} \}, \\
S_{k,6} & =\{\gamma : \hbox{ $\gamma$ has periodic expansion, period $(B_nB_s)^kB_nB_m^k$} \}, \\
S_{k,7} & =\{\gamma : \hbox{ $\gamma$ has periodic expansion, period $(B_nB_s)^k(B_n'B_s')^k $} \}, \\
S_{k,8} & =\{\gamma : \hbox{ $\gamma$ has periodic expansion, period $(B_nB_s)^k B_m'$} \}, \\
S_{k,9} & =\{\gamma : \hbox{ $\gamma$ has periodic expansion, period $(B_m)^k B_n'E_3'B_s'$} \}, 
\endalign $$
and write
$$ \delta_{\square}= M^*(\alpha,\gamma),\;\;\gamma \in S_{\square}. $$
In all cases except $a=3, b=4,5,6$ there are only finitely many values above the first
limit point. In those cases (due to small partial quotients) the behaviour is atypical
and we require additional sets of $\gamma$:

For $(a,b)=(3,4)$ define 
$$\align   S_{-8} & =\{\gamma : \hbox{ $\gamma$ has periodic expansion, period $F_0B_0$ or ${F_2B_2'} $ } \},\;\;\; \\
 S_{k,10} &  =\{\gamma : \hbox{ $\gamma$ has periodic expansion, period $F_2 (F_2B_2')^k $ or $F_2 (F_0B_0)^k $ } \},\;\;\; k\geq 1. \endalign  $$
For $(a,b)=(3,5)$ define blocks
$$G=\left(\frac{1}{2}(b-2-1),(a-1),\frac{1}{2}(b-2-1)\right), $$
$$ H=\left(\frac{1}{2}(a-2+1),\frac{1}{2}(b-2-3),(a-1),\frac{1}{2}(b-2-3),\frac{1}{2}(a-2+1)\right), $$
$$ H'=\left(\frac{1}{2}(a-2-1),\frac{1}{2}(b-2+1),(a-1),\frac{1}{2}(b-2+1),\frac{1}{2}(a-2-1)\right), $$
and sets
$$  S_{-9}  =\{\gamma : \hbox{ $\gamma$ has periodic expansion, period $HGH'G$ } \},\;\;\; $$
$$ S_{k,11}  =\{\gamma : \hbox{ $\gamma$ has periodic expansion, period $(H'GHG)^kHG$ } \},\;\;\; k\geq 0. $$
For $(a,b)=(3,6)$ define  
$$ S_{k,12}  =\{\gamma : 
\hbox{ $\gamma$ has periodic expansion, period $F_2 {B_2'}^{k+1}$ or 
$F_0B_2^kB_0$ } \},\;\;\; k\geq 0. $$

\proclaim{Theorem 2}
Suppose that $\alpha=[0;\overline{a,b}]$ with $a<b$ and $a\geq 3$ odd.
Then
$$ \rho^* (\alpha)=\cases \delta_{0}, & \hbox{ if $2\leq r\leq a-1$, $a\geq 5$ or $a=3$, $m\geq 2$, } \\
\delta_{-1}, & \hbox{ if $r=a+1$, $a\geq 5$ or $a=3$, $m\geq 1$}, \\
\delta_{-2}, & \hbox{ if $a+3\leq r\leq 2a$, } \\
 \delta_{-6}, & \hbox{ if $(a,b)=(3,4)$, } \\
\delta_{-7}, & \hbox{ if $(a,b)=(3,5)$. }
\endcases $$

The remaining values $\rho_2(\alpha),...$ down to the first limit point are
as follows. In all cases $\delta_{k,i}\nearrow \delta_{\infty,i}$
as $k\rightarrow \infty$ for the appropriate $i=1,...,9$.
The values of these $\delta_{k,i}$ and $\delta_{i}$ are given
in Lemma 2 below. 

\noindent
If $r\geq a+3$ then
$$ \rho_2^* (\alpha)=\cases \delta_{\infty,1}, & \hbox{ if $m\geq 1,$} \\
\delta_{\infty,2}, & \hbox{ if $m=0$, $a\geq 5$. } \endcases $$
If $r=a+1$ then 
$$ \rho_2^* (\alpha) =\cases \delta_{\infty,3}, & \hbox{ if $m\geq 1$, } \\
 \delta_{\infty,4}, &  \hbox{ if $m=0$, $a\geq 5$. } \endcases $$
If $4\leq r\leq a-1$ then
$$\align  \rho_2^* (\alpha) &  =\cases \delta_{-5},  &  \hbox{ if $b=a+4\geq 17$}, \\
  \delta_{\infty,5}, & \hbox{ otherwise, } \endcases \\
  \rho_{3}^*(\alpha)  & =\delta_{\infty,5}, \;\;\; \hbox{ if $b=a+4\geq 17$}.\endalign  $$
If $r=2$ then
$$ \align \rho_2^* (\alpha) & =\cases \delta_{-3},  & \hbox{ if $m\geq 2$, }\\
\delta_{-5}, & \hbox{ if $m=1$, $a\geq 9$,} \endcases \\
 \rho_3^* (\alpha) & =\cases \delta_{\infty,6}, & \hbox{ if $m\geq 3$,} \\
\delta_{\infty,7}, & \hbox{ if $m=2$,} \\
  \delta_{\infty,9}, & \hbox{ if $m=1$, $a\geq 9$.}  \endcases \endalign $$
For the remaining few cases:

\noindent
If  $(a,b)=(5,7)$ or $(7,9)$ then
$$\align \rho_2^* (\alpha) & =\delta_{-3},  \\
 \rho_3^* (\alpha) & =\delta_{-4}, \\
 \rho_4^* (\alpha) & =\delta_{\infty,8}.  \endalign $$

\noindent
If $(a,b)=(3,4)$, then the values above $\delta_{\infty,10}$ 
are $\delta_{-6}$, $\delta_{-8}$ and $\delta_{k,10}\searrow \delta_{\infty,10}$, $k\geq 0$. 

\noindent
If $(a,b)=(3,5)$ then the values above $\delta_{\infty,11}$ 
are $\delta_{-7}$, $\delta_{-9}$ and  $\delta_{k,11}\searrow \delta_{\infty,11}$, $k\geq 0$.   

\noindent
If $(a,b)=(3,6)$ then the values above $\delta_{\infty,12}$ 
are $\delta_{-2}$, $\delta_{-6}$ and  $\delta_{k,12}\searrow \delta_{\infty,12}$, $k\geq 0$.

\endproclaim 

We give the proof of Theorem 2 in $\S5$ below.

\head   The spectrum when $a=2$, $b\geq 5$. \endhead

We exclude here $b=3,4$ since the spectrum for $b=3$ is the same as for the period one $[0;3,\overline{4}]^-=1-[0;\overline{2,3}]^-$ considered in [6]  and the case $b=4,$ $[0;\overline{2,4}]=2-\sqrt{2}$ has been considered already by
Varnivides [\Var].

For $2\leq t\leq b-4$, with $t$ of the same parity as $b$, define
$$ S_{0,t} : =\left\{ \gamma \; :\; t_i 
\hbox{ periodic,  period $a,-t$ or $a,t-4$}\right\}. $$
For $b$ odd define
$$ \align S_{-2} & : =\left\{ \gamma \; :\; t_i 
\hbox{ period $(a,-3,a,-1)$}\right\}, \\
S_{-1} & : =\left\{ \gamma \; :\; t_i 
\hbox{ period $(a,-1)(a,-3,a,-1)(a,-1)$ or $(a,-3)(a,-1,a,-3)(a,-3)$}\right\}, \\
 S_{2k+1} & :=\left\{ \gamma \; :\; t_i 
\hbox{ period $(a,-1,0,-1)(a,-3,a,-1)^k$ 
  or $(a,-1,0,-1)(a,-1,a,-3)^k$}\right\} ,\;\; k\geq 0, \\
  S_{2k} & :=\left\{ \gamma \; :\; t_i 
\hbox{ period $(a,-1)(a,-3,a,-1)^k$ or $(a,-3)(a,-1,a,-3)^k$}\right\},\;\; k\geq 1, \endalign $$
and for $b$ even
$$ \align S_{-1} & : =\left\{ \gamma \; :\; t_i 
\hbox{ period $0,0$}\right\}, \\
 S_{2k+1} & :=\left\{ \gamma \; :\; t_i 
\hbox{ period $(a,-2,0,0)(a,-2)^k$ or $(a,0,0,-2)(a,-2)^k$}\right\} ,\;\; k\geq 0, \\
  S_{2k} & :=\left\{ \gamma \; :\; t_i 
\hbox{ period $(a,-4)(a,-2)^k$ or $(a,0)(a,-2)^k$}\right\} ,\;\; k\geq 1, \endalign $$
and set
$$ \delta_{k}:=M^*(\alpha,\gamma),\;\;\; \gamma \in S_{k} . $$

\proclaim{ Theorem 3}
Suppose that $\alpha =[0;\overline{2,b}]^-$ with $b\geq 5$.

Then for $b$ even
$$\rho^* (\alpha)=\cases \delta_{0,2}=\eta,  & \hbox{ if $b$ is even,}\\
\delta_{0,3}= \frac{2}{3}\eta, & \hbox{ if $b=5,$}\\
 \delta_{0,3}=\eta \left(1- \frac{1}{b^2-b}\right), & \hbox{ if $b\geq 7$ is odd.}\endcases $$
The remaining values of $M(\alpha,\gamma)$ greater than
$$ \delta_{\infty}=\eta \left( (1-\beta)^2 -\cases \beta^2, & \hbox{ if $b$ is even} \\  \frac{1}{4}\beta^4,  & \hbox{ if $b$ is odd.} \endcases \right), $$
are $\delta_{-1}$ ($b$ even), $\delta_{-2},\delta_{-1}$ ($b$ odd),
the  $\delta_{2k}$, $\delta_{2k+1}\searrow \delta_{\infty}$, and  when $b\geq 8$ the $\delta_{0,m}$ with $4 \leq m\leq 2+[\sqrt{2b-4}]$ 
($m$ of the same parity as $b$). 
For these values  $M^*(\alpha,\gamma)=\delta_{j}$ iff $\gamma$
is in $S_{j}$.

\endproclaim

\head  3. The value of  the points in the spectrum. \endhead

The following Lemmas give the values, $\delta_{i}$, described in the above Theorems.

\proclaim{Lemma 1}
Suppose that $a\geq 4$ is even.

When $b$ is odd
$$ \delta_0=\left(1-\beta -\frac{1}{b}\right)\left(1-\eta +\frac{\eta}{b}
\right), $$
$$ \delta_{-1}=\left(1-\beta -\frac{\beta (1-D)}{1+D^2}\right)
\left(1-\eta -\frac{D(1-D)}{1+D^2}\right), $$
$$\delta_{-2}=\cases \left(1-\beta +\left(\frac{3\beta -2D}{1-D}\right)\right)\left(1-\eta -\left(\frac{2\eta -3D}{1-D}\right)\right), & \hbox{ if $b\geq \max\{2a-5,3a/2\},$} \\
\left(1-\beta -\left(\frac{3\beta -2D}{1-D}\right)\right)\left(1-\eta +\left(\frac{2\eta -3D}{1-D}\right)\right), & \hbox{ if $b\leq \min\{a+5,3a/2\},$} \\
\left(1+\beta -\left(\frac{3\beta -2D}{1-D}\right)\right)\left(1+\eta -\left(\frac{2\eta -3D}{1-D}\right)\right), & \hbox{ otherwise,} \endcases 
$$
and
$$\delta_{k}=\left( 1-\beta -\frac{\beta (1-D)(1+2D^2)}{1+D^2}-\beta D^3 w_{k}\right)\left(1-\eta -\frac{D(1-D)}{1+D^2}+w_{k}\right) $$
where $w_{k}:=2D^{4k+1}(1-D)/((1+D^2)(1-D^{4k+2})). $

When $b$ is even,
$$\delta_{k,1}=\left(1-\eta +\frac{2D^2}{1+D}(1- w_{k,1})\right)
\left(1-\beta -\frac{2\beta }{1+D}(1- w_{k,1})\right)$$
where $w_{k,1}:=D^{2k}(1-D)/(1-D^{2k+1})$, $k\geq 0$.

For $b=2a-2$, $a\geq 8$,
$$\delta_{k,2}=\left(1-3\eta+\frac{2D(2-\eta)}{1-D}-w_{k,2}\right)\left(1+\beta -\frac{2\beta D}{1-D}(1-\eta+D)+\beta D w_{k,2}\right) $$
where $w_{k,2}:=2D^{k+1}(1+D)(1-\eta+D)/((1-D)(1+D^{k+2})).$

For $b=2a-4$, $a\geq 10$,
$$\align \delta_{k,3} =  & \left(1-\eta-\frac{2\eta D}{1-D} +\frac{2D}{1-D^2}(1+2D)+w_{k,3}\right) \\
   & \hskip3ex \left(1-3\beta +\frac{2D}{1-D} -\frac{2\beta D}{1-D^2}(2+D)-\beta D w_{k,3}\right) \endalign $$
with $w_{k,3}:=2D^{2k+1}/((1+D)(1-D^{2k+1}))$.

For $a+6\leq b\leq 2a-6$
$$ \delta_{k,4}=\left( 1-\eta +\frac{2D(1-\eta)}{1-D} 
+ w_{k,4}\right)\left(1-3\beta +\frac{2D (1-\beta)}{1-D}-\beta w_{k,4}\right) $$
where $w_{k,4}:=2D^{k+1}/(1-D^{k+1}).$

For $k=0$ or $b\leq 2a-6$ or $a=6,b=8$, or $k=1$ and $b=2a-4$
$$ \delta_{k,5}=   \left( 1-\eta +\frac{2D(1-\eta)}{1-D} +\eta w_{k,5}\right)  \left(1-3\beta +\frac{2D (1-\beta)}{1-D}-w_{k,5}\right)  $$
where $w_{k,5}:=2D^{k+1}(1-2\beta+D)/((1-D)(1+D^{k+1})).$

For $(a,b)=(8,12)$
$$\align \delta_{k,6}=  & \left( 1-3\eta +2D-2D^2+2\eta D^2-\frac{2D^3(1-\eta+D)}{1+D^2} -\eta D^2 w_{k,6}\right)  \\
  & \hskip3ex \left(1+\beta -\frac{2D (1-\beta+\beta D)}{1+D^2}+w_{k,6}\right) 
\endalign $$
where $w_{k,6}:=2D^{4k+4}(1-2\beta+D)(1-D^3)/((1+D^2)(1-D^{4k+6})).$

For $(a,b)=(6,10)$
$$\align \delta_{k,7}=  & \left( 1+\eta -2D+2D^2+\frac{2\eta D^3}{1-D}-\frac{2D^3(1+2D)}{1-D^2} -2\eta D^3w_{k,7}\right)  \\
   & \hskip3ex \left(1-5\beta +\frac{2D}{1-D}-2\beta D\frac{(1+2D)}{(1-D^2)}-w_{k,7} \right) \endalign  $$
where $w_{k,7}:=2D^{2k}(1-3\beta +D)/(1-D^{2k+2})$, and $\delta_{1,4}=\frac{703}{40}-\frac{703}{600}\sqrt{210}$.

For the remaining $\delta_{i}$ we have
$$ \delta_{-1}=\left(1-3\eta+\frac{2D(1-\eta)}{1-D}\right)\left(1-\beta +\frac{2\beta (1-\eta)}{1-D}\right) $$
and
$$\delta_{0,4}=\cases \left(1+3\beta-\frac{2D(1-2\beta)}{1-D}\right)\left(1-3\eta +\frac{2D(2-\eta)}{1-D}\right), & \hbox{ if $b\geq \max\{ 3a-6, 2a\}$,} \\
\left(1-5\beta+\frac{2D(1-2\beta)}{1-D}\right)\left(1+\eta -\frac{2D(2-\eta)}{1-D}\right), & \hbox{ if $b\leq\min\{ a+6, 2a-2\}$,} \\ 
\left(1-3\beta+\frac{2D(1-2\beta)}{1-D}\right)\left(1-\eta +\frac{2D(2-\eta)}{1-D}\right), & \hbox{ else,} 
 \endcases $$
and for $(a,b)=(4,6)$
$$ \delta_{-2}=\left(1-\eta -\frac{2D}{1-D}+\frac{2\eta D}{1-D^2}\right)\left(1-3\beta -\frac{2\beta D}{1-D}+\frac{2D}{1-D^2}\right). $$
\endproclaim

\proclaim{Lemma 2}
Suppose that $a$ is odd.

Set
$$ v:=\frac{m\beta -D}{1-D},\;\;\;\; D:=\eta\beta. $$
Then
$$ \delta_{0}=(1-2\eta+\eta v)(1-\beta+v), \tag 3.1 $$
$$ \delta_{-1}=\cases  \left(1-\eta v-\frac{2D^2}{1-D^2}\right)\left(1-3\beta -v -\frac{2\beta D^2}{1-D^2}\right), & \hbox{ if $r\leq a+1$, } \\
 \left(1-2\eta +\eta v+\frac{2D}{1-D^2}\right)\left( 1-\beta +v +\frac{2\beta D}{1-D^2}\right), & \hbox{ if $r\geq a+1$,} \endcases \tag 3.2  $$
$$ \delta_{-2}=\left(1-\eta v-\frac{2D}{1-D}\right)\left(1-3\beta -v -\frac{2\beta D}{1-D}\right), \tag 3.3 $$
$$\delta_{-3}= \left(1-2\eta +\eta v +\frac{2\eta}{b}\right)\left(1-3\beta +v +\frac{2D}{b}\right), \tag 3.4 $$
$$\delta_{-4}= \left(1-2\eta +\frac{\eta (m+\eta)}{b}\right)\left(1-\beta -\frac{(m-\eta)}{b}\right). \tag 3.5 $$
If $m=1$ then
$$\delta_{-5}=\cases \left(1-4\eta +\frac{3D (1-\eta)}{1-D}\right)
\left(1+2\beta-\frac{3D(1-\beta)}{1-D}\right), & \hbox{ if $r\geq a-7$}, \\
\left( 1-2\eta +\frac{3D(1-\eta)}{1-D}\right) \left( 1-2\beta +\frac{3D(1-\beta)}{1-D}\right), & \hbox{ if $r\leq a-9$.} \endcases \tag3.6 $$
For $b$ even
$$ \delta_{-6}=\eta $$
and for $b$ odd
$$\delta_{-7} =\eta \left(1-\left(\frac{\beta}{1-D}\right)^2\right).$$
For $r\geq a+3$
$$\delta_{k,1}=\left( 1-2\eta +\eta v +\frac{2D}{1-D}-\frac{2D^{k+1}}{1-D^{k+1}}\right)\left(1-\beta+v +\frac{2\beta D}{1-D} -\frac{2\beta D^{k+1}}{1-D^{k+1}}\right). \tag 3.7 $$
For $m=0$, $k\geq 1$, $r\geq a+3,$
$$\delta_{k,2}=\left(1-2\eta +\frac{D(2-\eta)}{1-D} -\eta \ve_{k,2}\right)\left(1-\beta +\frac{D(1-2\beta)}{1-D}+D\ve_{k,2}\right). \tag 3.8  $$
For $r\leq a+1$ and $b\geq 6$
$$ \delta_{k,3}=\left(1-\eta v -\frac{2D}{1-D^2}-\eta \ve_{k,3}\right)\left(1-3\beta -v -\frac{2\beta D^2}{1-D^2}-\ve_{k,3}\right). \tag 3.9 $$ 
For $b=a+1\geq 6$
$$\delta_{k,4}=\left(1-\frac{\eta D}{1-D}+\frac{2D^2}{1-D^2}+\eta D \ve_{k,4}\right)\left(1-3\beta +\frac{D}{1-D}-\frac{2\beta D^2}{1-D^2} -\ve_{k,4}\right). \tag 3.10$$
 For $r\leq a-1$
$$ \delta_{k,5}=\left(1-\eta v-\frac{2D^{k+1}}{1-D^{k+1}}\right)\left(1-3\beta -v -\frac{2\beta D^{k+1}}{1-D^{k+1}}\right). \tag 3.11 $$  
For $r=2$ and $b\geq 7$
$$ \delta_{k,6}=\left(1-\eta v-\frac{2D^{k+1}(1+D^{2k+1})}{(1+D)(1-D^{3k+1})}\right)\left(1-3\beta -v +\frac{2D}{b}-\frac{2\beta D^{2k+1}(1+D^k)}{(1+D)(1-D^{3k+1})}\right). \tag 3.12 $$
For $b=2a+2$
$$\delta_{k,7}=\left(1-\frac{\eta D}{1-D} +\frac{4D^2}{1-D^2} +\eta D\ve_{k,7}\right)\left(1-\beta +\frac{D}{1-D} -\frac{4\beta}{1-D^2} -\ve_{k,7}\right). \tag 3.13 $$
For $b=a+2\geq 7$
$$\align  \delta_{k,8}=  & \left(1+\frac{D(1+D-4D^2)}{1-D^2} -\frac{\eta D (1-2D)}{1-D} -\eta D^2 \ve_{k,8}\right)  \\
  & \hskip3ex \left(1-4\beta +\frac{D}{1-D} +\frac{\beta D(1-3D)}{1-D^2}-\ve_{k,8}\right). \tag 3.14 \endalign   $$
and $b=a+2\geq 11$
$$\align \delta_{k,9}=   & \left(1-2\eta +3D-3\eta D+3D^2-\eta D^2-\frac{\eta D^2(\beta -D)}{1-D}-\eta D^2\ve_{k,9}\right)\\
 & \hskip3ex \left(1-2\beta +\frac{D(1-\beta)}{1-D}-\ve_{k,9}\right). \tag 3.15 \endalign $$
Here the $\ve_{k,i}\rightarrow 0$ are given explicitly by
$$  \ve_{k,2}:=\frac{2D^{k}\left(\beta (1+D)-D\right)}{(1-D)(1+D^{k+1})}, \;\;\; \ve_{k,3}:=\frac{2\beta D^{2k+1}}{(1+D)(1-D^{2k+1})}, $$ 
$$ \ve_{k,4}:=\frac{2D^{2k}(1-2/b)}{(1-D)(1+D^{2k})},\;\;\; \ve_{k,7}:=\frac{2D^{2k}(1-4/b)}{(1-D)(1+D^{2k})}, $$
$$\;\;\;
\ve_{k,8}:=\frac{2D^{2k}(1-2/b)}{1-D^{2k+1}},\;\;\; \ve_{k,9}:=\frac{2D^{k+1}}{1-D^{k+3}}(1-2\beta+2D-2\beta D+D^2).$$

For $(a,b)=(3,4)$
$$ \delta_{-8}=\eta \left(1 -\left(\frac{\beta (1-\eta+D)}{1-D^2}\right)^2\right), $$
$$ \delta_{k,10}=\eta \left( 1-\frac{\beta (1-\eta +D)}{1-D^2}+\ve_{k,10}\right)\left(1+\frac{\beta D (1-\eta +D)}{1-D^2}-D\ve_{k,10}\right) $$
where $\ve_{k,10}=\beta D^{2k}(1-\eta+D)/ (1+D)(1-D^{2k+1})$. 

For $(a,b)=(3,5)$
$$\delta_{-9}=\eta \left(1-\left(\frac{2\beta -D+D^3-2\beta D^3}{1+D^4}\right)^2\right), $$
$$\align  \delta_{k,11}= & \eta \left( 1-2\beta +D -D^3\frac{(1-2\beta -2\beta D+D^2)}{1+D^4}(1-\ve_{k,11})\right)\\
  &  \left( 1+2\beta -D -D^3\frac{(1-2\beta -2\beta D+D^2)}{1+D^4}(1+D^4 \ve_{k,11})\right) \endalign $$
where $\ve_{k,11}=2 D^{8k}/(1-D^{8k+4}).$
\endproclaim

For $(a,b)=(3,6)$
$$\delta_{k,12}=\eta \left(1-\left(\frac{\beta (1-\eta +D)(1-D^{k+1})}{(1-D)(1-D^{k+2})}\right)^2\right). $$

\proclaim{Lemma 3}
Suppose that $a=2$.
For $2\leq t\leq b$
$$\delta_{0,t} = \cases \eta \left(1-\left(\frac{(t-2)\beta}{1-D}\right)^2\right), & \hbox{ if $t\leq b-\sqrt{2b-4},$} \\
\eta \left( \left(2-2\beta -\frac{(t-2)\beta}{1-D}\right)^2 -1\right), & \hbox{ if $t> b-\sqrt{2b-4}$.} \endcases  $$
For $b$ even
$$ \align \delta_{-1} & =\eta (1-\beta)^2,  \\
\delta_{2k} & =\eta \left(1- 2 \beta -\frac{2\beta D^{k+1}}{1-D^{k+1}}\right)\left(1+\frac{2\beta D^k}{1-D^{k+1}}\right),  \\
\delta_{2k+1} & =\eta \left( \left(1-  \beta\right)^2 -\beta^2 \left(\frac{1-D^{k+1}}{1-D^{k+2}}\right)^2\right). \endalign $$
For $b$ odd
$$ \align \delta_{-2} & =\eta \left(1-\beta +\frac{\beta D}{1+D}\right)^2, \\
\delta_{-1} & =\eta \left( \left(1-\beta +\frac{\beta (1-D)D^3}{1-D^4}\right)^2-\left(\frac{\beta (1+D)D}{1-D^4}\right)^2\right),   \\
\delta_{2k} & =\eta \left(1-  \beta -\frac{1}{2}\beta^2 -\frac{2\beta D^{2k+2}}{(1+D)(1-D^{2k+1})}\right)\left(1-\beta +\frac{1}{2}\beta^2 +\frac{2\beta D^{2k}}{(1+D)(1-D^{2k+1})}\right),  \\
\delta_{2k+1} & =\eta \left( \left(1-  \beta\right)^2 -\frac{1}{4}\beta^4 \left(\frac{1-D^{2k}}{1-D^{2k+2}}\right)^2\right). \endalign $$

\endproclaim

\head  Proof of Lemma 1 (even $a\geq 4$) \endhead

Notice that (because of the forwards backwards symmetry, 
and the rough bound $s_4^*(2i-1)\geq (1+\eta -D)(1+\beta )>1$ when $d_{2i-1}^+>0$) in all cases 
it is enough to
check
$$u_i:=\min \{s_1^*(2i-1),s_2^*(2i-1),s_3^*(2i-1)\}, \;\;\; \;\;\;d_{2i-1}^+\geq 0, $$
for each $i$ in the period of $\gamma$ or its negative, where
$s_1^*(2i-1)<s_2^*(2i-1)$ iff $(d_{2i-1}^+-\beta)(1+d_{2i-1}^-)<\eta$,
with $s_1^*(2i-1)<s_3^*(2i-1)$ iff $d_{2i-1}^-(1-\beta)+d_{2i-1}^+(1-\eta)<0$ and
$s_2^*(2i-1)<s_3^*(2i-1)$ iff $(d_{2i-1}^-+\eta)(1-d_{2i-1}^+)+\beta <0$.

For period $C_3$ we have $d_{2i-1}^+=(3\beta -2D)/(1-D)$,
$d_{2i-1}^-\sym (-2\eta +3D)/(1-D)$ and one tediously checks which of the three
functions gives $u_i$.

When $(b_{2i-1},b_{2i})=A_0$, $A_1$, $A_2$, or $C_2$ the bound $d_{2i-1}^+\leq 2\beta$ 
plainly gives $s_1^*(2i-1)<s_2^*(2i-1)$.

For blocks $C_2$ the bounds $d_{2i-1}^-\leq -2\eta+4D$, $d_{2i-1}^+\leq 2\beta$
give $s_1^*(2i-1)<s_3^*(2i-1)$ and hence $u_i=s_1^*(2i-1)$.

For blocks $A_1$ or $A_2$ the bounds $d_{2i-1}^-\geq -2D$, $d_{2i-1}^+\geq \beta -D$ give $u_i=s_3^*(2i-1)$.

Hence the values of $\delta_{0,1}$, $\delta_{0,5}$ and $\delta_{-1}$ with $b$ even, and $\delta_0$ when $b$ is odd.

Observe that if a block $A_1$ is preceded by an $A_1'$ then $d_{2i-1}^+\leq \beta +\beta D$, $d_{2i-1}^-\leq -D+2D^2$ gives $s_3^*(2i-1)\geq (1-2\beta -\beta D)(1-\eta +D-2D^2)>(1-2\beta)(1-\eta)$ greater than the value claimed for
$\delta_{-1}$, $b$ odd and the  $\delta_{k}$, $k\geq 1$.
Hence in these cases we need only check the $A_1$ preceded by an $A_1$.
For $\delta_{-1}$ these have $d_{2i-1}^+=\beta (1-D)/(1+D^2)$, $d_{2i-1}^-\sym
D (1-D)/(1+D^2)$, and $s_3^*(2i-1)$ gives $\delta_{-1}$. 
For blocks $A_1'A_1(A_1'A_1'A_1A_1)^k$ observe that taking the negative
essentially reverses the sequence (interchanging $d_{2i-2}^-$ and $d_{2i-1}^+$)
and we need just check $s_3:=(1-\eta -\eta d)(1-\beta -d')$ where
$d=\min\{ d_{2i-2}^-,d_{2i-1}^+\}$, $d'=\max\{ d_{2i-2}^-,d_{2i-1}^+\}$
for $t_{2i}=1=t_{2i-2}$. Plainly the largest $d,d'$ 
both occur when the $A_1$ is followed by $A_1'A_1A_1'$ and 
$$d'\sym \frac{\beta (1-D)}{1-D^{4k+2}}\left(1+\frac{D^2(1-D^{4k})}{1+D^2}\right),\;\;\; d\sym \frac{\beta  (1-D)}{1-D^{4k+2}}
\left(\frac{(1-D^{4k})}{1+D^2}-D^{4k}\right),$$
with $s_3^*(2i-1)$ giving the value claimed for $\delta_{k}$, $k\geq 1$.

For the blocks $A_0$ in $S_{k,1}$, $k\geq 1$, the bounds $d_{2i-1}^+\leq 2\beta D$, $d_{2i-1}^-\leq -2D+2D^2$ give $u_i=s_1^*(2i-1) \sym \left(1-\eta -\frac{2D}{1+D}(1-w_{k,1})\right)
\left(1-\beta +\frac{2\beta D}{1+D} (1-w_{k,1})\right), $
while cutting at the $A_2$ the minimum value of the  $u_i=s_3^*(2i-1)$ occurs 
when $d^+= 2\beta (1-w_{k,1})/(1+D)$, $d^-\sym -2D^2(1-w_{k,1})/(1+D)$ 
(i.e. when $d^-$ and $d^+$ are both largest) with this
 giving the value claimed for $\delta_{k,1}$.

For the $\delta_{k,5}$ observe that cutting at an $A_2$ gives
$$d^-\sym \left(-2D(1-\eta)\frac{(1-D^k)}{(1-D)}-2D^{k+1}\right)/(1+D^{k+1}),$$
$$ d^+=\left(2\beta -2D(1-\beta)\frac{(1-D^k)}{(1-D)}\right)/(1+D^{k+1})$$
and $u_i=s_3^*(2i-1)$ asymptotically gives the value claimed for $\delta_{k,5}$ with these increasing to $\delta_{\infty,5}=\delta_{\infty,4}$. Cutting at a $C_2$ gives
$d^-\geq (-2\eta+2D)/(1-D)$ and $d^+\geq (2\beta-2D)/(1-D)$ and 
$$ \align u_i=s_1^*(2i-1) & \geq \left(1-\eta -\frac{2\eta}{1-D}(1-\beta)\right)\left(1-\beta +\frac{2\beta}{1-D}(1-\eta)\right) \\
  & =\delta_{\infty,5}+\frac{2D}{(1-D)^2}(2a-5-b+7\eta-4D-\eta D+D^2) >\delta_{\infty,5}\endalign $$
for $b\leq 2a-6$ or $a=6, b=8$.
For $k=1$ and $b=2a-4$ one checks that these
$$ u_i=s_1^*(2i-1)\sym \left(1-3\eta +2D\frac{(1-D+\eta D)}{1+D^2}\right)
\left(1+\beta -2\beta D\frac{(1-\eta+D)}{1+D^2}\right)$$
are still greater than the value claimed for $\delta_{1,5}$;
$$ u_i=s_3^*(2i-1)\sym \left(1-\eta +2D\frac{(1-\eta +D)}{1+D^2}\right)\left(1-3\beta +2 D\frac{(1-\beta+\beta D)}{1+D^2}\right).$$


Cutting at a $C_4$ the bounds 
$$\frac{(4\beta-2D)}{(1-D)}\geq d^+\geq (4\beta -2D),\;\;\; \frac{(-2\eta +4D)}{(1-D)}\geq d^-\geq -2\eta +2D-2\eta D, $$
give $s_1^*(2i-1)\leq s_3^*(2i-1)$ iff  $b\geq 2a$,
give $s_2^*(2i-1)<s_3^*(2i-1)$ if  $b\geq a+8$ (or $b=a+6$ when $C_4$
is not preceded by another $C_4$, using  $d^-\leq -2\eta+2D$) and
$s_2^*(2i-1)> s_3^*(2i-1)$ if  $b\leq a+4$ (or $b=a+6$ when $C_4$
is  preceded by another $C_4$, using $d^-\geq -2\eta+4D-2\eta D$),
and give $s_1^*(2i-1)<s_2^*(2i-1)$ if $b\geq \max\{3a-6, 2a\}$
and $s_2^*(2i-1)<s_1^*(2i-1)$ if $b\leq 3a-8$. Hence the expressions for $\delta_{0,4}$.

For the $S_{k,2}$ observe that cutting at a $C_2$ produces
$$ d^+=\left(2\beta -2\beta D +2D^2(1-2\beta)\frac{(1-D^k)}{(1-D)}+2D^{k+2}\right)/(1+D^{k+2}) $$
$$d^-\sym \left(-2\eta +2D(2-\eta)\frac{(1-D^k)}{(1-D)}+(2-2D)D^{1+k}\right)/(1+D^{k+2}) $$
and $u_i=s_1^*(2i-1)$ asymptotically gives the value claimed for 
$\delta_{k,2} $ with these
increasing to $\delta_{\infty,2}\leq 1-3\eta +\beta +D-2\eta D+2\beta D+10D^2.$
Cutting at an $A_2$ the bounds 
$d^-\leq -2D+2\eta  D$, $d^+\leq 2\beta -2D+4\beta D$, give
$$ \align u_i=s_3^*(2i-1) \geq & (1-\eta +2D-2\eta D)(1-3\beta +2D-4\beta D) \\
  \geq & 1-\eta-3\beta +7D-4\eta D-10\beta D +10D^2 \\
 \geq &  \delta_{\infty,2}+2D(b-2a+3-2\eta-4\beta) >\delta_{\infty,2}\endalign $$
for $b=2a-2$, $a\geq 6$. Cutting at a $C_4$ the bounds $d^+\leq 4\beta -2D+4\beta D$, 
$-2\eta +4d\geq d^-\geq -2\eta +2D-2D^2$, give
$$\align  s_2^*(2i-1) & \geq (1-\eta +2D-2D^2)(1-3\beta+2D-4\beta D)\\
  & \geq 1-3\beta-\eta+7D-2\eta D-10\beta D 
   \geq \delta_{\infty,2}+2D(b-2a+3 -\eta -4\beta-5D) \endalign $$
(and $s_3^*(2i-1)>\delta_{\infty,2}$ for $(a,b)=(8,14)$). Hence the value of $\delta_{k,2}$ for $b=2a-2,$ $a\geq 8.$

For the $\delta_{k,3}$ observe that for the $C_4$ the smallest $d^-$
and largest $d^+$ both occur when
$$ d^+=\frac{4\beta -2D}{1-D}-\frac{2\beta D^2 (1-D^{2k})}{(1-D^2)(1-D^{2k+1})},\;\;\; d^-\sym \frac{-2\eta +4D}{1-D}-\frac{2D(1-D^{2k})}{(1-D^2)(1-D^{2k+1})}, $$
so that when $a+6\leq b<2a$ we have $u_i=s_2^*(2i-1)\leq (1-\eta+2D)(1-3\beta+2D)$ 
giving the value claimed for
$\delta_{k,3}$ (checking numerically that for $(a,b)=(10,16)$ one
has $u_i=s_3^*(2i-1)\geq (1+\eta-4D+2\eta D-2D^2)(1-5\beta +2D-2\beta D)$
 for a $C_4$ preceded by a $C_4$, greater than the claimed value $\delta_{k,3}$). Cutting at a $C_2$ the rough bounds
$d^-\geq -2\eta+4D-2\eta D$, $d^+\geq 2\beta -2D$ give
$ u_i=s_1^*(2i-1)\geq (1-3\eta+4D-2\eta D)(1+\beta -2D)$,
which is larger when $b=2a-4$.

For the $\delta_{k,4}$ cutting at the $C_4$ gives
$$ d^+=\frac{2\beta-2D}{1-D}+\frac{2\beta}{1-D^{k+1}},\;\;\; d^-\sym \frac{-2\eta+2D}{1-D}+\frac{2D^{k+1}}{1-D^{k+1}} $$
and when $a+6\leq b <2a$
$$u_i= s_2^*(2i-1)\sym \left(1-\eta+\frac{2D(1-\eta)}{1-D}+w_{k,4}\right)\left(1-3\beta +\frac{2D(1-\beta)}{1-D} -\beta w_{k,4}\right)$$
gives the value claimed. Cutting at a $C_2$
$$\align  u_i=s_1^*(2i-1) & \geq \left(1-3\eta +\frac{2D-2\eta D}{1-D}\right)\left( 1+\beta -\frac{(2D-2\beta D)}{1-D}\right) \\
  & \geq \left(1-\eta +\frac{4D-2\eta D}{1-D}\right)\left(1-3\beta +\frac{2D-4\beta D}{1-D}\right) \endalign $$
(the value claimed when $k=0$) for $b\leq 2a-6$. 


For $(a,b)=(8,12)$ and $\delta_{k,6}$ cutting at the first $C_2$ in a block $C_2C_2$ gives
$$ d^-\sym \left(-2\eta +2D +2\eta D^2 -2D^2-2D^3 \frac{(1-\eta+D)(1-D^{4k+4})}{1+D^2}\right)/(1-D^{4k+6}),\;\;\;$$
$$ d^+\sym \left(2\beta \frac{(1-\eta +D)(1-D^{4k+4})}{1+D^2}+ 2D^{4k+4}(1-\beta +\beta D -D^2) \right)/(1-D^{4k+6})$$
and $u_i=s_1^*(2i-1)$  asymptotically gives the value claimed for $\delta_{k,6}$ 
(the largest value corresponding to $k=0$). For the remaining $C_2$ and $A_2$
we have the even larger
$$ u_i=s_1^*(2i-1)\geq (1-3\eta+2D-2\eta D)(1+\beta -2\beta D) $$
and
$$ u_i=s_3^*(2i-1)\geq (1-\eta +2D-2\eta D)(1-3\beta+2D-2\beta D). $$
 
For $(a,b)=(6,10)$ cutting at the $C_4$ in blocks $A_2C_4$
gives $s_3^*(2i-1)\leq \delta_{\infty,7}$ asymptotically equal to the value claimed
for $\delta_{k,7}$. For the remaining $C_4$
$$ u_i=s_3^*(2i-1)\geq (1-5\beta+2D-2\beta D)(1+\eta-2D+2\eta D-4D^2)>\delta_{\infty,7}$$
and for the $A_2$ and $C_2$ respectively we have
$$ u_i=s_3^*(2i-1)\geq (1-\eta +2D-2D^2)(1-3\beta+2D-4\beta D+2D^2-2\beta D^2)>\delta_{\infty,7},$$
$$ u_i=s_1^*(2i-1) \geq \left(1+\beta-2D+4\beta D-2D^2 \right)\left(1-3\eta +4D-2\eta D \right)>\delta_{\infty,7}. $$

Finally when $(a,b)=(4,6)$, in the case of  $\delta_{-2}$, cutting at an $A_2$ we have $u_i=s_3^*(2i-1)$ 
asymptotically giving the value claimed, while cutting
at the $C_2$ gives 
$$ u_i=s_1^*(2i-1)\sym \left(1-3\eta +\frac{2D}{1-D}-\frac{2\eta D^2}{1-D^2}\right)\left(1+\beta  +\frac{2\beta D}{1-D}-\frac{2D^2}{1-D^2}\right). \;\;\; \blacksquare$$


\head  Proof of Lemma 2 (odd $a$) \endhead

We first exclude the $S_{k}$ containing $t_{2i-1}=a$.
Hence it is enough to check the minimum of the $s_j^*(2i),s_{j}^*(2i+1)$
for $j=1,...,4$ for all $i$ in the period of $\pm \gamma$ with $t_{2i+1}<0$.
Since the $d_{2i}^-,d_{2i+1}^+\geq -2\beta -D$ we can ignore
the $s_4^*(2i+1),s_2^*(2i)\geq (1+\eta)(1-\beta -D)>(1-\eta)(1-\beta)$ and
 need only check
$$ v_i:=\min\{ s_1,s_3, s_2 \}, \hskip4ex t_{2i+1}<0, $$
where, writing $d:=\min\{d_{2i+1}^+,d_{2i}^-\}$, $d':=\max\{d_{2i+1}^+,d_{2i}^-\}, $
$$\align  s_1 & :=\min\{s_1^*(2i+1),s_1^*(2i)\}=(1-\eta +t_{2i+1}\eta +\eta d')(1-\beta +d), \\
  s_3 & :=\min\{s_3^*(2i+1),s_3^*(2i)\}=(1-\eta -t_{2i+1}\eta -\eta d)(1-\beta -d'),\\
 s_2 & :=\min\{s_2^*(2i+1),s_4^*(2i)\}=(1+\eta +t_{2i+1}\eta +\eta d)(1+\beta -d').\endalign $$
Now if $t_{2i+1}=-1$ then $v_i=\min\{s_1,s_3\}$ (if $d'>0$ then $d\geq -2\beta-D$ gives $-\eta d(1-d')<\beta $ and $s_3<s_2$,
while if $d'<0$ then plainly $s_1<s_2$). 
For $\gamma$ in $S_{-4}$ we have $d=-(m\beta -D)/(1+D)$, $d'=(m\beta+D)/(1+D)$
and $v_i=s_1$ gives the value of $\delta_{-4}$.

We first deal with the $\gamma$ having $t_{2i+1}=-1$ for all $i$.
Writing
$$ d=\frac{m\beta -D}{1-D}+2\beta w, \;\;\; d'=\frac{m\beta -D}{1-D}+2\beta w', $$
for $t_{2i+1}=-1$ we have 
$$ s_1<s_3 \hbox{ iff }   w(1-\eta -D)+w'(1+\eta -D) <\frac{r}{a}. $$
In particular 
$$ \{t_{2i},t_{2i+2}\}=\{m,m\} \Rightarrow   v_i=s_1$$
giving (3.1) for $\delta_0$. Likewise
$$ \{t_{2i},t_{2i+2}\}=\{n,n\} \Rightarrow   v_i=\cases s_1, & \hbox{ if $r=2a$,}\\ s_3, & \hbox{  if $r\leq 2a-2$}, \endcases $$
(with asymptotic equality when $r=2a$ and $t_{2i}=n$ for all $i$) and (3.3) is clear.
Similarly
$$ \{t_{2i},t_{2i+2}\}=\{n,m\} \Rightarrow   v_i=\cases s_1, & \hbox{ if $r\geq a+3$,}\\ s_3, & \hbox{  if $r\leq a-1$}, \endcases $$
where when $r=a+1$ we have $v_i=s_3$ if $w'\geq 1/(1-D^2)$, $w\geq D/(1-D^2)$
and  $v_i=s_1$ if  $w'\leq 1/(1-D^2)$, $w\leq D/(1-D^2)$
(with asymptotic equality when the $t_{2i}$ have period $n,m$)
giving the value (3.2) of $\delta_{-1}$.
For the $\gamma$ in $S_{-3}$ we have $w'=1/(1+D)$, $w=-1/(1+D)$ and $s_1$
gives $\delta_{-3}$.

Suppose $r\geq a+3$ and $\gamma$ is in $S_{k,1}$. 
If $\{t_{2i},t_{2i+2}\}=\{m,n\}$ then 
$d'=v+2\beta/(1-D)-2\beta D^k/(1-D^{k+1})$ and
$d=v+2\beta/(1-D)-2\beta/(1-D^{k+1}) $ and
$$ \align v_i= s_1 & 
= \left(1-2\eta+\eta v+\frac{2D}{1-D}-\frac{2D^{k+1}}{1-D^{k+1}}\right)
\left(1-\beta +v+\frac{2\beta D}{1-D} -\frac{2\beta D^{k+1}}{1-D^{k+1}}\right)
 \\
  & <\delta_{\infty,1}:=\left(1-2\eta+\eta v+\frac{2D}{1-D}\right)
\left(1-\beta +v+\frac{2\beta D}{1-D}\right). \endalign $$
For the remaining $\{t_{2i},t_{2i+2}\}=\{n,n\}$ we have 
$v+2\beta<d,d'<v +2\beta/(1-D)$ and
$$\align s_1 & \geq (1-2\eta+\eta v+2D)
(1-\beta+v+2\beta ) \\
  & > \delta_{\infty,1}+\frac{2\beta}{1-D}(1-2D)(1-2\eta)-\frac{2D^2}{1-D}
(1+\eta)>\delta_{\infty,1},\endalign $$
while
$$ \align
s_3 & \geq 
\left(1-\eta v-\frac{2D}{1-D}\right)\left(1-\beta -v 
-\frac{2\beta}{1-D}\right) \\
  & =\delta_{\infty,1}+2\beta 
\left( \frac{r-(a+3)}{a}+\eta \left(1-\frac{3\beta}{a}+\frac{\beta(m+2)-D}{1-D}\right)\right) >\delta_{\infty,1},\endalign $$
and (3.7) holds.

Suppose that $r\leq a-1$ and $\gamma$ is in $S_{k,5}$. 
If $\{t_{2i},t_{2i+2}\}=\{m,n\}$ we have 
$d'=v+2\beta/(1-D^{k+1})$, 
$d=v+2\beta D^k/(1-D^{k+1})$ and
$$ \align 
v_i=s_3 & =
\left(1-\eta v-\frac{2D^{k+1}}{1-D^{k+1}}\right)
\left(1-3\beta -v-\frac{2\beta D^{k+1}}{1-D^{k+1}} \right) \\
 & <\delta_{\infty,5}:=(1-\eta v)(1-3\beta -v).\endalign $$
For the remaining $\{t_{2i},t_{2i+2}\}=\{m,m\}$ we have 
$d,d'\geq v$ and
$$\align  v_i=s_1 & \geq 
(1-2\eta+\eta v)(1-\beta+v) \\
  & = \delta_{\infty,5}+2\beta \left( \frac{a-1-r}{a} +\frac{1}{a}-
\eta v\right) >\delta_{\infty,5}, \endalign
 $$
giving (3.11).

Suppose that $r\leq a+1$ and $\gamma$ is in $S_{k,3}$. If $\{t_{2i},t_{2i+2}\}=\{n,n\}$ we have 
$d,d'=v+2\beta/(1-D^2)+\ve_{k,3}$ and
$$ \align v_i=s_3 & =
\left(1-\eta v-\frac{2D}{1-D^2}-\eta \ve_{k,3}\right)
\left(1-3\beta -v -\frac{2\beta D^2}{1-D^2}-\ve_{k,3}\right) \\
 & < \delta_{\infty,3}:=\left(1-\eta v-\frac{2D}{1-D^2}\right)
\left(1-3\beta -v -\frac{2\beta D^2}{1-D^2}\right). 
\endalign  $$
For the other $\{t_{2i},t_{2i+2}\}=\{n,m\}$ we have 
$d'\leq v +2\beta/(1-D)$ and
$d \leq v+2\beta D/(1-D)$ and for $b\neq 4,5$
$$\align v_i=s_3 &  \geq 
\left(1-\eta v-\frac{2D^2}{1-D}\right)
\left(1-\beta-v-\frac{2\beta}{1-D}\right) \\
 & \geq \delta_{\infty,3}+\frac{2D}{1-D^2}(1-(4+m)\beta)>\delta_{\infty,3},
\endalign $$
giving (3.9) for $b\neq 4,5$.

Suppose that $r=2$ and that $\gamma$ is in $S_{k,6}$. If 
$\{t_{2i},t_{2i+2}\}=\{n,m\}$ then
$$ d'=v+\frac{2\beta (1+D^{2k+1})}{(1+D)(1-D^{3k+1})},\;\;\; 
d=v+\frac{2\beta D^{k}(1+D^{2k+1})}{(1+D)(1-D^{3k+1})}$$
and $v_i=s_3$ equals 
$$ \align  &\left(1-\eta v 
-\frac{2D^{k+1}(1+D^{2k+1})}{(1+D)(1-D^{3k+1})}\right)
\left(1-\beta -v -\frac{2}{b}-
\frac{2\beta D^{2k+1} (1+D^k)}{(1+D)(1-D^{3k+1})}\right) \\
  & < \delta_{\infty,6}:=(1-\eta v)\left(1-\beta -v 
-\frac{2}{b}\right).\endalign $$
If $\{t_{2i},t_{2i+2}\}=\{m,m\}$ we have 
$d,d'>v$ and
$$\align v_i =s_1 & \geq 
\left(1-2\eta +\eta v\right)(1-\beta +v) \\
 & =\delta_{\infty,6} +
\frac{2}{b}(1-2\eta -\eta v)
\endalign $$
If $\{t_{2i},t_{2i+2}\}=\{s,n\}$ then $w>-1/(1+D)$, $w'>1/(1+D)$
and 
$d' \leq v+2\beta /(1-D)$, $d\leq v-2\beta +2\beta D/(1-D)$ give
$$\align
v_i=s_3 & \geq 
\left(1-\eta v+\frac{2D(1-2D)}{1-D}\right)\left(1-\beta -v 
-\frac{2\beta}{1+D}-\frac{4\beta D}{1-D}\right) \\
 & > \delta_{\infty,6}-\frac{4\beta D}{1-D}(1-\eta v+2D)+
\frac{2D(1-2D)}{1-D}(1-3\beta -v) \\
& > \delta_{\infty,6}+\frac{2D}{1-D}(1-(5+m)\beta -D) >
\delta_{\infty,6}\endalign $$
for $b>5$, giving (3.12).

For the remaining cases we have to deal with non constant $t_{2i+1}$.

Suppose that $\gamma$ is in $S_{k,2}$, $k\geq 1$, with $m=0$ and $b=r\geq a+3$.
If $\{t_{2i},t_{2i+2}\}=\{n,-m\}$ then
$d'=(2\beta-D)/(1-D)-\ve_{k,2}$, $d=(D-2\beta D)/(1-D)+D\ve_{k,2}$
and $w'\leq 1/(1-D)$, $w\leq \eta$ gives
$$ \align
v_i=s_1 & =\left(1-2\eta +\frac{(2D-\eta D)}{(1-D)}-\eta \ve_{k,2}\right)\left(1-\beta +\frac{D-2\beta D}{1-D} +D\ve_{k,2}\right) \\
  & < \delta_{\infty,2}:=\left(1-2\eta +\frac{2D-\eta D}{1-D}\right)\left(1-\beta +\frac{D-2\beta D}{1-D}\right). \endalign $$
If $\{t_{2i},t_{2i+2}\}=\{n,n\}$ then $2\beta -D\leq d,d'\leq (2\beta -D)/(1-D)$
and
$$ \align s_3 & \geq \left(1-\frac{(2D-\eta D)}{1-D}\right)\left(1-\beta -\frac{(2\beta -D)}{1-D}\right) \\
  & = \delta_{\infty,2}
+\frac{2D}{(1-D)^2}(b-a-3+(7\beta -\eta-2D-\beta D+D^2)) \endalign $$
and
$$ \align
s_1 & \geq (1-2\eta+2D-\eta D)(1-\beta+2\beta-D) \\
 & > \delta_{\infty,2} -\frac{2D^2}{1-D}(1+\beta)+(1-2\eta)\frac{(2\beta -2D)}{1-D} >\delta_{\infty,2},\endalign $$
and (3.8) holds.

Suppose $b=a+1$ and $\gamma$ is in $S_{k,4}$. Then at all points 
$w'\geq 1/(1-D^2)$, $w\geq D/(1-D^2)$ and $v_i=s_3$.
If $t_{2i+1}=-1$ and $t_{2i-1}$ or $t_{2i+3}=1$  then $d'=2\beta/(1-D^2)- D/(1-D)+\ve_{k,4}$, $d=D/(1-D)-2\beta D/(1-D^2)-D\ve_{k,4}$ and
$$ \align
s_3 & =\left(1-\frac{\eta D}{1-D}+\frac{2D^2}{1-D^2}+\eta D \ve_{k,4}\right)\left(1-\beta -\frac{2\beta}{1-D^2}+\frac{D}{1-D} -\ve_{k,4}\right)\\
 & < \delta_{\infty,4}:=\left(1-\frac{\eta D}{1-D}+\frac{2D^2}{1-D^2}\right)\left(1-\beta -\frac{2\beta}{1-D^2}+\frac{D}{1-D}\right).\endalign $$
If $t_{2i+1}=t_{2i-1}=t_{2i+3}=-1$ then $d'\leq 2\beta -D+D^2$, $d\leq -D+2\beta D$ and
$$\align s_3 & \geq (1+\eta D-2D^2)(1-3\beta +D-D^2) \\
  & \geq \delta_{\infty,4} +(2\eta D-4D^2)(1-3\beta+D)-2D^2(1+D) \\
  & \geq \delta_{\infty,4} +2\eta D(1-6\beta +D) >\delta_{\infty,4} \endalign $$
for $b\neq 4$ and (3.10) holds.

Suppose that $b=2a+2$ and that $\gamma$ is in $S_{k,7}$. Then in all cases
$w'>1/(1+D),w>-1/(1+D)$ and $v_i=s_3$.
If $t_{2i+1}=-1$ and $t_{2i-1}$ or $t_{2i+3}=1$ then $d'=4\beta/(1-D^2)-
D/(1-D)+\ve_{k,7}$ and $d=D/(1-D)-4\beta D/(1-D^2)-D\ve_{k,7}$ giving
$$\align 
s_3 & = \left(1-\frac{\eta D}{1-D}+\frac{4D^2}{1-D^2}+\eta D \ve_{k,7}\right)\left( 1-\beta +\frac{D}{1-D} -\frac{4\beta}{1-D^2}-\ve_{k,7}\right) \\
  & > \delta_{\infty,7}:=\left(1-\frac{\eta D}{1-D} +\frac{4D^2}{1-D^2} \right)\left(1-\beta +\frac{D}{1-D}-\frac{4\beta}{1-D^2}\right). \endalign $$
If $t_{2i+1}=t_{2i-1}=t_{2i+3}=-1$ 
then $d'\leq 4\beta -D+D^2$ and $d\leq -D+4\beta D$ and
$$ \align
s_3 & > (1+\eta D-4D^2)(1-5\beta +D-D^2) \\
& > \delta_{\infty,7} -2D^2(1+D)+(2\eta D-8D^2)(1-5\beta +D) \\
& > \delta_{\infty,7} +2\eta D(1-10\beta +D +17\beta^2) >\delta_{\infty,7}
\endalign $$
and (3.13) holds.

Next suppose that $b=a+2\geq 7$ and that $\gamma$ is in $S_{k,8}$. If $t_{2i+1}=-1$ and $t_{2i-1}$ or $t_{2i+3}=-1$ then $w'\geq 1/(1+D)$, $w\geq -1/(1+D)$ and $v_i=s_3$.
If $t_{2i+1}=-1$ and  $\{t_{2i-1},t_{2i+3}\}=\{-1,1\}$ then
$$ d'=\frac{\beta (3-D)}{1-D^2} -\frac{D}{1-D} +\ve_{k,8} ,\;\;\; d=-\beta -\frac{D}{1-D} +2D -\frac{\beta D(1-3D)}{1-D^2} +D^2\ve_{k,8},$$
and $v_i=s_3$ equals
$$\align 
&  \left(1+D +\frac{D^2 (1-3D)}{1-D^2} +\frac{\eta D^2}{1-D}-\eta D-\eta D^2 \ve_{k,8}\right)\left(1-\beta -\frac{\beta (3-D)}{1-D^2}+\frac{D}{1-D}-\ve_{k,8}\right) \\
 & < \delta_{k,8}:=\left(1+D+\frac{D^2 (1-3D)}{1-D^2} +\frac{\eta D^2}{1-D}-\eta D\right)\left(1-\beta -\frac{\beta (3-D)}{1-D^2}+\frac{D}{1-D}\right). \endalign $$
When $\{t_{2i-1},t_{2i+3}\}=\{-1,-1\}$ we have $d'\leq 3\beta -D -\beta D+D^2$
and $d \leq -\beta -D +3\beta D -D^2$ and
$$\align 
s_3 & \geq (1+D+\eta D -3D^2+\eta D^2)(1-4\beta +D +\beta D-D^2) \\
  & \geq \delta_{\infty,8} -2D^2(1+D+\eta D)+(2\eta D-4D^2)(1-4\beta +D) \\
  & > \delta_{\infty,8} +2\eta D (1-7\beta +D) \geq \delta_{\infty,8}
\endalign $$
for $b\geq 7$.  If $\{t_{2i-1},t_{2i+3}\}=\{1,1\}$ then $w,w'\leq 2D$ and $d,d'\geq \beta $ give $v_i=s_1\geq (1-2\eta +D)$
while
$\delta_{\infty,8} < (1+D)(1-4\beta +D+\beta D+D^2)< 1-4\beta +2D<1-2\eta +D$,
so (3.14) holds.

Suppose now that we also allow $t_{2i+1}=-3$ with $m=1$. For these points in $S_{-5}$ and $S_{k,9}$ we have $d,d'\leq 3\beta-D$ and $s_3-s_2=2\eta (2-d)(1-d')-2\beta >0$
with $s_1-s_2=2(d(1-3\eta -D+\eta d)-\eta -\beta +3D)$. Hence if $m=1$ and $\gamma$ is in $S_{-5}$ we have $d,d'\rightarrow \frac{3\beta-3D}{1-D}$ and
 writing
$$ s_1-s_2=\frac{2 D}{(1-D)^2}\left(a-r-9+3\beta +12\eta -9D +2\beta D-\eta D\right),$$
clearly  $v_i=s_1$
if $r\geq a-7$ and  $v_i=s_2$
if $r\leq a-9$, giving (3.6).
Suppose that $\gamma$ is in $S_{k,9}$ with $m=1$, $r=2$ and $a\geq 9$.
If $t_{2i+1}=-1$ and $\{t_{2i},t_{2i+2}\}=\{3,-1\}$ then
$$ d'=3\beta -3D+3\beta D-D^2-\frac{D^2 (\beta -D)}{1-D}-\ve_{k,9}D^2,\;\;\;
d=\frac{-\beta +D}{1-D}-\ve_{k,9}, $$
and $w'\leq 1-\eta +D$, $w\leq -1+\eta$ gives $v_i=s_1$ equal to
$$\align & \left(1-2\eta +3D-3\eta D+3D^2-\eta D^2 -\frac{\eta D^2 (\beta -D)}{1-D}-\eta D^2 \ve_{k,9}\right)\left(1-2\beta +\frac{D(1-\beta)}{1-D}-\ve_{k,9}\right) \\
  & <  \delta_{k,9} :=\left(1-2\eta +3D-3\eta D+3D^2-\eta D^2 -\frac{ D^3 (1-\eta )}{1-D}\right)\left(1-2\beta +\frac{D(1-\beta)}{1-D}\right). \endalign $$
For $t_{2i+1}=-3$ we have $d,d'\geq 3\beta -D-2\beta D$ 
giving $s_1-s_2\geq 2D(a-9+4\eta-7\beta +d^2/\beta )>0$ for $a\geq 9$
 and
$d,d'\leq 3\beta -D$ gives 
$$v_i=s_2\geq \left(1-2\eta +3D-\eta D\right)\left(1-2\beta +D\right) >\delta_{\infty,9}.  $$
Finally, when $t_{2i+1}=-1$, $t_{2i}=t_{2i+2}=1$ we have $w,w'<\eta$ and $d,d'\geq \beta -D$ 
give
$$ v_i=s_1\geq (1-2\eta +D-\eta D)(1-D)>1-2\eta >\delta_{\infty,9}. $$

For the configurations $a=3$, $b=4,5,6$ we need also to deal with $t_{2i+1}=a$.
For $t_{2i+1}=a$ we have
$$ s_4:=\min\{s_2^*(2i),s_4^*(2i+1)\}=\eta (1-\beta -d')(1+\beta +d). $$
From the bounds $d,d'\geq  -2\beta +D-\beta D-\beta D^2$ when $b=4$, 
$d,d'\geq -3\beta +D-D^3$ for $b=5$ and  $d,d'\geq -2\beta $ for $b=6$
we have $(2-3\beta+\beta^2)+3d+d^2+(d'-d)(1+d)>0$ and
$s_1:=\min\{s_1^*(2i),s_1^*(2i+1)\}=\eta (5-\beta+d')(1-\beta+d)>s_4$
while $s_2^*(2i+1)=(1+4\eta+\eta d_{2i}^-)(1+\beta-d_{2i+1}^+)>\eta$.
Noting that switching from $\gamma$ to $-\gamma$ replaces
$d_{2i}^-$ and $d_{2i+1}^+$ by $-d_{2i}^- -2\beta$ and $-d_{2i+1}^+-2\beta$
when $t_{2i+1}=a$ (and by  $-d_{2i}^-$ and $-d_{2i+1}^+$ as before otherwise)
it is readily seen that we need just evaluate $s_4$ for the  $t_{2i+1}=a$ with 
$t_{2i}<0$, and $\min\{s_1,s_2,s_3\}$ for the $t_{2i+1}=-1$ (in  both $\gamma$ and $-\gamma$).
For  $\gamma$ in $S_{-6}$ and $S_{-7}$ we $d,d'\sym -\beta$ and $d,d'\sym -\beta-\beta /(1-D)$  and $\delta_{-6}$ and $\delta_{-7}$  are immediate.
For $(a,b)=(3,4)$ and $\gamma$ in $S_{-8}$ cutting at the $t_{2i+1}=a$ gives
$$ d,d'\sym \frac{-2\beta +D-2\beta D+aD^2}{1-D^2}$$
with $s_4$ giving the value claimed for $\delta_{-8}$ while for $\gamma$ in $S_{k,10}$ cutting at the $t_{2i+1}=a$ with $\{t_{2i-1},t_{2i+3}\}=\{a,1\}$ gives
$$ d'=\left(-2\beta +aD +\frac{(-2\beta D+D^2-2\beta D^2+aD^3)(1-D^{2k})}{1-D^2}\right)/(1-D^{2k+1}), $$
$$ d=\left( \frac{(-2\beta +D-2\beta D+aD^2)(1-D^{2k})}{1-D^2}-2\beta D^{2k}+aD^{2k+1} \right)/(1-D^{2k+1}), $$
and the value claimed for $\delta_{k,10}$, the remaining $t_{2i+1}=a$ 
producing larger $d$ and smaller $d'$.
It is readily checked that cutting at a $t_{2i+1}=-1$ 
we have $\beta \leq d,d'\leq \beta (1+D-\eta D)/(1-D^2)$ 
and $\min \{ s_1,s_2,s_3\}\geq \eta $ does not affect the value.

Similarly for  $(a,b)=(3,5)$ when cutting at the $t_{2i+1}=a$, $t_{2i}=t_{2i+2}=-3$,
for  $\gamma$ in $S_{-9}$ we have
$d,d' \sym -3\beta +D -D^3(1-2\beta-2\beta D+D^2)/(1+D^4)$ 
while  for $S_{k,11}$ the smallest $d$ and largest $d'$ for $\gamma$ occur when
$t_{2i+7}$ or $t_{2i-5}=1$ and
$d\sym -3\beta+D -D^3(1-2\beta-2\beta D+D^2)(1-\ve_{k,11})/(1+D^4)$
and $d'\sym -3\beta+D +D^3(1-2\beta-2\beta D+D^2)(1+D^4\ve_{k,11})/(1+D^4)$.
These give the values claimed for $\delta_{-9}$ and $\delta_{k,11}$.
When $t_{2i+1}=a$ and $t_{2i}=t_{2i+2}=-1$ we have $d'\leq -\beta +D-2\beta D-2\beta D^2+D^3+D^4$, $d\geq -\beta -D+2\beta D+2\beta D^2-D^3-D^4$ and $s_4> \eta \left(1-\left(\frac{2\beta-D-D^3+2\beta D^3}{1-D^4}\right)^2\right)$ (the value claimed for $\delta_{0,11}$ and the largest of those claimed). For the $t_{2i+1}=-1$ we have $-D^2\leq d\leq D^2$ and $2\beta +\beta D\leq d'\leq 2\beta +2\beta D$ and $\min\{s_1,s_2,s_3\}$ is also larger. 
When $k\geq 1$  the $t_{2i+1}=a$, $t_{2i}=t_{2i+2}=-3$ in $-\gamma$ have
$d\geq -3\beta +D-D^3+2\beta D^3+2\beta D^4-D^5-D^7$ and 
$d'\leq -3\beta +D-D^3+2\beta D^3+2\beta D^4-D^5+D^7$ and $s_4$ is certainly
greater than the value claimed for $\delta_{1,11}$.

For $(a,b)=(3,6)$ and $\gamma$ in $S_{k,12}$ cutting at the $t_{2i+1}=a$
gives 
$$ d,d' \sym -\beta -\frac{\beta (1-\eta +D)(1-D^{k+1})}{(1-D)(1-D^{k+2})} $$
and the value claimed for $\delta_{k,12}$, while the remaining $t_{2i+1}=-1$ have $\beta \leq d,d'\leq 2\beta$ and
$\min\{s_1,s_2,s_3\}\geq \eta$.
 $\;\;\; \blacksquare$

\head  Proof of Lemma 3 ($a=2$)  \endhead   

We assume that the $a_{2i-1}=a=2$ and the $a_{2i}=b\geq 5$. Using the forwards-backwards symmetry in the expansions of the $\gamma$ in question it is straightforward to see that it is enough
to check
$$\align 
s_1^*(2i-1)  & =  (1-\eta +d_{2i-1}^-)(1-\beta +d_{2i-1}^+), \\
s_2^*(2i-1)  & =  (1+\eta +d_{2i-1}^-)(1+\beta -d_{2i-1}^+), \\
s_4^*(2i-1)  & =  (1+\eta -d_{2i-1}^-)(1+\beta +d_{2i-1}^+), \\
\endalign $$
and the corresponding functions, $s_1^*(2i-1)',s_2^*(2i-1)',s_4^*(2i-1)'$
say, for  $-\gamma$.

Now if $t_{2i-1}=0$ then  going from $\gamma$ to $1-\alpha-\gamma$ gives
$({d_{2i-1}^-},{d_{2i-1}^+})'=(-d_{2i-1}^-, -d_{2i-1}^+),$
and $s_1^*(2i-1)',s_2^*(2i-1)',s_4^*(2i-1)'$ merely add the function
$$ s_3^*(2i-1)=  (1-\eta -d_{2i-1}^-)(1-\beta -d_{2i-1}^+). $$
Writing $s_1^*(2i-1),s_3^*(2i-1)=\eta (1-\beta \pm d_{2i-2}^-)(1-\beta\pm d_{2i-1}^+)$
clearly
$$ v_i:=\min\{ s_1^*(2i-1),s_3^*(2i-1)\} =\eta (1-\beta -v')(1-\beta +v) $$
where 
$$v':=\max\{-d_{2i-2}^-,d_{2i-1}^+  \},\hskip5ex v:=\min\{-d_{2i-2}^-,d_{2i-1}^+\}.$$
Now for the $\gamma$ in question we have $|d_{2i-1}^+|\leq \beta,|d_{2i-1}^-|\leq 
D$, giving $s_2^*(2i-1),s_4^*(2i-1)>1$ (certainly larger than $\delta_{k}$), so it is enough just to
evaluate $v_i$ when $t_{2i-1}=0$.

Now if $t_{2i-1}=a$ then
$$ s_4^*(2i-1)=\eta (1-\beta -d_{2i-2}^-)(1+\beta +d_{2i-1}^+) $$
while going from $\gamma$ to $1-\alpha-\gamma$ gives
$({d_{2i-2}^-},{d_{2i-1}^+})'=(-2\beta -d_{2i-2}^-,-2\beta -d_{2i-1}^+)$,
$$ s_4^*(2i-1)'=\eta (1+\beta +d_{2i-2}^-)(1-\beta -d_{2i-1}^+), $$
and
$$ u_i:=\min\{s_4^*(2i-1),s_4^*(2i-1)'\}=\eta (1+\beta +d)(1-\beta-d'), $$
where
$$ d':=\max\{ d_{2i-2}^-, d_{2i-1}^+    \},\hskip5ex  d:=\min\{ d_{2i-2}^-, d_{2i-1}^+    \}. $$
Suppose that the $t_i$ in the expansion of $\gamma$ have period $a,-t$.
Then  $d_{2i-1}^-\sym (a\eta -tD)/(1-D)$, $d_{2i-2}^-,d_{2i-1}^+\sym (-t\beta +aD)/(1-D)$ and $d_{2i-1}^+(1+\eta)\leq d_{2i-1}^-(1+\beta)$
giving  $s_4^*(2i-1)\leq s_2^*(2i-1)$. Now if $t\geq 2$ then it is easily seen that
${d_{2i-1}^+}',{d_{2i-1}^-}'\geq {d_{2i-1}^+},{d_{2i-1}^-}$ so that
$$ \min\{s_1^*(2i-1),s_1^*(2i-1)'\}=s_1^*(2i-1)\sym \eta \left( \left(2-2\beta -\frac{(t-2)\beta}{1-D}\right)^2-1\right) $$
while 
$ u_i \sym \eta \left(1-\left( \frac{(t-2)\beta}{1-D}\right)^2\right) $
and $\delta_{0,t}$ will be  the minimum of these two expressions.
Note that the second expression is smaller than the first  when $t\leq b-\sqrt{2b-4}$
and greater than $\delta_{\infty}$ when $(t-2)\leq \sqrt{2b-4}$
(for the finite number of $b$ and $t$ with $(t-2)\leq \sqrt{2b-4}$
but $t> b-\sqrt{2b-4}$ one can check numerically that  $\delta_{0,t}\geq
\delta_{\infty}$ except for $(b,t)=(6,4)$ or $(7,5)$).

For the remaining $\gamma$ if $t_{2i-1}=a$ we have $d_{2i-1}^+<aD$, $d_{2i-1}^-\geq 1-\eta$ so $s_4^*(2i-1)<s_2^*(2i-1)$. Also the rough bounds
$d_{2i-2}^-,d_{2i-1}^+\geq (-4\beta +aD)/(1-D)$
when $b\geq 8$ and $d_{2i-2}^-,d_{2i-1}^+\geq (-3\beta +aD)/(1-D)$
when $b=7$, and $d_{2i-2}^-, d_{2i-1}^+\geq (-4\beta +aD-2\beta D+aD^2)/(1-D^2)$, $(-2\beta +aD-4\beta D+aD^2)/(1-D^2)$ for $\gamma$ in $S_{2k}$
and $d_{2i-2}^-,d_{2i-1}^+\geq -2\beta$ for  $\gamma$ in $S_{2k+1}$
 when $b=6$,
give $\beta <\eta (1+d_{2i-2}^-)(1+d_{2i-1}^+)=(1+d_{2i-1}^+)(d_{2i-1}^- -\eta)$ and hence $s_4^*(2i-1)\leq s_1^*(2i-1)$ when $b\geq 6$.
When $b=5$ we have $d_{2i-1}^-\geq (a\eta -3D)/(1-D)$, $d_{2i-1}^+\geq (-3\beta +aD)/(1-D)$ and $s_1^*(2i-1)\geq 2\eta /3$ (larger than the values claimed for the
$\delta_{k}$) and so can be similarly ignored. 
Hence for the $t_{2i-1}=a$ we need only consider $u_i$.

For period $(a,-4)(a,-2)^k$ the smallest $d$ and largest $d'$ 
occur simultaneously when $t_{2i}=-4$ (or $t_{2i-2}=-4$) and
$$\align d' & =d_{2i-2}^-\sym \frac{-2\beta +aD}{1-D}-\frac{2\beta D^k}{1-D^{k+1}}=-\beta -\frac{2\beta D^k}{1-D^{k+1}} \\
d & =d_{2i-1}^+ = \frac{-2\beta +aD}{1-D}-\frac{2\beta }{1-D^{k+1}}=-3\beta -\frac{2\beta D^{k+1}}{1-D^{k+1}} \endalign $$
with $u_i$ asymptotically giving the value claimed for $\delta_{2k}$, $b$ even.

Similarly for period $(a,-3)(a,-3,a,-1)^k$ the smallest $d$ and largest
$d'$ again occur simultaneously when $...,t_{2i-2},t_{2i-1}=...-3,a,-3,a$
and 
$$\align d & =d_{2i-2}^- \sym \frac{aD}{1-D} +\left( -3\beta +(-3\beta D-\beta D^2)\frac{(1-D^{2k})}{1-D^2}\right)/(1-D^{2k+1}) \\
 & =-2\beta -\frac{1}{2}\beta^2 -\frac{2\beta D^{2k+2}}{(1+D)(1-D^{2k+1})} \\
d' & =d_{2i-1}^+ =\frac{aD}{1-D} +\left( (-\beta -3\beta D)\frac{(1-D^{2k})}{1-D^2} -3\beta D^{2k}\right)/(1-D^{2k+1}) \\
 &  = -\frac{1}{2}\beta^2 -\frac{2\beta D^{2k}}{(1+D)(1-D^{2k+1})} \endalign $$
with this $u_i$  asymptotically giving the value claimed for $\delta_{2k}$, $b$ odd.

For period $(a,-3,a,-1)$ we clearly have $d'\sym (-\beta +aD-3\beta D+aD^2)/(1-D^2)$ and $d\sym (-3\beta +aD-\beta D+aD^2)/(1-D^2)$,
giving $\delta_{-2}$. Likewise for $\gamma$ in $S_{-1}$, $b$ odd, the largest and smallest $d'$, $d$
occur when cutting the sequence $...,t_{2i-1}|t_{2i},...$
at the point $-3,a,-3,a,-3,a|-1,a$ and
$$ \align d' & =d_{2i-1}^+=\frac{-\beta -3\beta D -3\beta D^2-3\beta D^3}{1-D^4}+\frac{aD}{1-D}=-\frac{\beta D (1+D+D^2-D^3)}{1-D^4}\\
d & =d_{2i-2}^-\sym \frac{-3\beta -3\beta D -3\beta D^2-\beta D^3}{1-D^4}+\frac{aD}{1-D}=-2\beta -\frac{\beta D (1+D-D^2+D^3)}{1-D^4}, \endalign $$
with $u_i$ giving the stated value of $\delta_{-1}$.

Now if $\gamma$ has period  $(a,0,0,-2)(a,-2)^k$ then cutting at a $t_{2i-1}=a$
gives
$$ \align
d' & \leq  \frac{-2\beta D +aD^2}{1-D^2}=-\frac{1}{2}\beta^2, \\
d & \geq  \frac{-2\beta +aD^2}{1-D^2}=-2\beta +\frac{1}{2}\beta^2,
\endalign $$
and the $u_i\geq \eta (1-\beta +\frac{1}{2}\beta^2)^2$, greater than the value claimed for $\delta_{2k+1}$, $b$ even.

Now cutting at the $t_{2i-1}=0$ gives
$$\align
d_{2i-2}^- & \sym \frac{aD-2\beta D}{1-D} + \frac{(-aD^2+2\beta D^2)D^k}{1-D^{k+2}} =  \frac{\beta (1-D^{k+1})}{1-D^{k+2}} \\
d_{2i-1}^+ & =  \frac{-2\beta +a D}{1-D} + \frac{(2\beta D-aD^2)D^k}{1-D^{k+2}} = - \frac{\beta (1-D^{k+1})}{1-D^{k+2}} 
\endalign $$
and $v_i$ gives the valued claimed for  $\delta_{2k+1}$, $b$ even.

Finally for period $(a,-1,0,-1)(a,-1,a,-3)^k$ cutting at a $t_{2i-1}=a$ we have
$$\align 
d  & \geq \frac{-3\beta +aD-\beta D-\beta D^2+a D^3-\beta D^3+a D^4}{1-D^4} =-2\beta +\frac{\beta (1-D)D^3}{1-D^4} \\
d'  & \leq \frac{-\beta +aD-\beta D-\beta D^2+a D^3-3\beta D^3+a D^4}{1-D^4} 
=- \frac{\beta (1-D)D^3}{1-D^4} \endalign $$
and $u_i\geq \eta (1-\beta)^2$ is larger than the value 
claimed for $\delta_{2k+1}$.
Lastly cutting at the $t_{2i-1}=0$ we have
$$\align
d_{2i-1}^+ & = -\beta + \frac{aD -\beta D +a D^2 -3\beta D^2}{1-D^2} - \frac{(a D-2\beta D)D^{2k+1}}{1-D^{2k+2}} = \frac{1}{2}  \frac{\beta ^2 (1-D^{2k})}{1-D^{2k+2}} \\
d_{2i-2}^- & \sym  \frac{-\beta + aD -3\beta D +a D^2}{1-D^2} + \frac{(2\beta -a D)D^{2k+1}}{1-D^{2k+2}} = -\frac{1}{2}  \frac{\beta ^2 (1-D^{2k})}{1-D^{2k+2}} \endalign $$
and $v_i$ asymptotically gives the value claimed for $\delta_{2k+1}$, 
$b$ odd.
$\blacksquare$


\head 4. Proof of Theorem 1 (even $a\geq 4$) \endhead

Writing
$$ u=\cases 1-\eta -2\beta +D, & \hbox{ if $b$ is odd,} \\
1-\eta -3\beta +3D, & \hbox{ if $b$ is even,} \endcases $$
we note that $\delta_{\infty}\geq (1-2\beta)(1-\eta -D)>u$ when $b$ is odd
and $\delta_{\infty,1} \geq (1-\eta)(1-3\beta)=u $ when $b$ is even.

Since $u>\eta$ we can by (1.2) rule out $t_i=a_i$ infinitely often. In
particular we can assume  that
going from $\gamma$ to $1-\gamma$ merely replaces the $t_i$ 
with $-t_i$. Since $1-4\eta +D<u$ we can further assume 
from (1.3) that $|t_{2i-1}|=0,2$ for all $i$, and since $1-2\eta +D<u$
when $b$ is odd and $\eta >2\beta $ or when $b$ is even and $\eta \geq 3\beta -2D$,
that $t_{2i-1}=0$ when $b\geq 2a+1$ is odd
or $b\geq 3a-2$ is even. Likewise we can assume that $|t_{2i}|=1,3$ when
$b\leq 2a-1$ is odd and $|t_{2i}|=0,2$ or $4$ when $b\leq 3a+2$ is even
and $|t_{2i}|=0,2$ when $b=a+2$
(else $M^*(\alpha,\gamma)$ is less than $1-5\beta +D<u$ or $1-6\beta +D<u$
or $1-4\beta +D<u$ respectively).
From the bound 
$$\delta_{\infty,4}\geq (1-\eta+2D-2\eta D)(1-3\beta+2D-2\beta D)\geq 1-3\beta-\eta+7D-4\eta D -8\beta D $$
and numerically for $\delta_{\infty,6}$ we similarly obtain $|t_{2i}|=0,2$ for $b=a+4$, $a\geq 8$.

Now if $t_{2i+1}=0$ and $d:=\max\{|d_{2i+1}^+|,|d_{2i}^-|\}\geq 2\beta$
when $b$ is odd or $d\geq 3\beta$ when $b$ is even then
$$\min\{ s_1^*(2i+1), s_3^*(2i+1), s_1^*(2i), s_3^*(2i)\}\leq (1-\beta -d)(1-\eta +\eta d)  \tag 4.1 $$
is less than $u$.
Hence we can assume if $t_{2i+1}=0$ then $t_{2i},t_{2i+2}=\pm 1$ 
or $0,\pm 2$ as $b$ is odd or even.

If $t_{2i+1}=-2$ then (with at most a finite number of exceptions) 
$t_{2i},t_{2i-2}\geq 0$ else (since the $t_{2i}\leq 4$) the $s_1^*(2i)$ or $s_1^*(2i+1)$ give
$$ \align M^*(\alpha,\gamma) & \leq \left(1-3\eta +\frac{4D+2\eta D}{1-D}\right)\left(1-3\beta +\frac{2D+4\beta D}{1-D}\right) \\
   & \leq 1-3\eta-3\beta+15D-4\eta D<u \endalign $$
for $b$ even and
$$ \align  M^*(\alpha,\gamma) & \leq \left(1-3\eta +\frac{3D+2\eta D}{1-D}\right)\left(1-2\beta +\frac{2D+3\beta D}{1-D}\right) \\
  & \leq 1-3\eta -2\beta +11D-4\eta D <u \endalign $$
for $b$ odd  (in particular if $t_{2i}=3$ or $4$ we can assume that $t_{2i\pm 1}=-2$).

Suppose first that $b$ is odd. If 
$t_{2i+1}=-2$ and $t_{2i}\neq 3$ then $t_{2i}=1$
and
$$ s_1^*(2i)\leq \left(1-3\eta +\frac{3D-2\eta D}{1-D}\right)\left(1-\beta +\frac{\beta}{1-D}\right)< 1-3\eta +3D-2\eta D +2\beta D\leq u. $$
Hence for $b$ odd we can assume that the expansion of $\gamma$ either has period $C_3$ or $C_3'$ or consists solely of blocks $A_1'$, $A_1$.
We suppose the expansion consists only of blocks $A_1, A_1'$ and assume that $\gamma$ is not in $S_{0}$ so that $\gamma$ (or its negative) has infinitely
many blocks $t_{2i-1},...,t_{2i+3}=0,-1,0,-1,0$. We can rule out blocks
$A_1A_1A_1$ or $A_1'A_1'A_1'$ since $t_{2i+4}=-1$ gives
$$ s_1^*(2i+1)\leq \left(1-\eta -D+\frac{D^2}{1-D}\right)\left(1-2\beta -\beta D+\frac{\beta D^2}{1-D}\right)< (1-2\beta)(1-\eta -D), $$
less than $\delta_{\infty}$, and we can rule out $t_{2i+4},t_{2i+5},...$ (or $t_{2i-2},t_{2i-3},...$)
of the forms $1,0,-1,0,-1,...$ or $1,0,-1,0,1,0,-1,...$ else
$$ d_{2i+1}^+\leq -\beta +\beta D-\beta D^2+\frac{\beta D^3(1-D+D^2)}{1-D^3} =-\frac{\beta (1-D)(1+2D^2)}{1+D^2}-\frac{2\beta D^4 (1-D-D^3)}{(1+D^2)(1-D^3)}, $$
$$ d_{2i+1}^-\leq \frac{-D+D^2+D^3}{1-D^3}=-\frac{D(1-D)}{1+D^2}+\frac{2D^5}{(1+D^2)(1-D^3)}, $$
and $s_1^*(2i+1)\leq \delta_{\infty}-\frac{2\beta D^4}{(1+D^2)(1-D^3)} ( (1-D-D^3)(1-\eta -D)-\eta )<\delta_{\infty}$.
Hence either $\gamma$ has period $A_1A_1A_1'A_1'$ or $\gamma$ (or its negative)
has infinitely many blocks $t_{2i-1},...,t_{2i+11}=0,-1,0,-1,0,1,0,-1,0,1,0,1,0$
and writing $d_{2i+1}^-=-D(1-D)/(1+D^2)+w^-$ and 
$d_{2i+8}^+=D(1-D)/(1+D^2)-w^+$ we have
$$\align s_1^*(2i+1) & = \left(1-\eta -\frac{D(1-D)}{1+D^2}+w^-\right)\left(1-2\beta +\beta D -\frac{\beta D^2 (1-D)}{1+D^2}-\beta D^3 w^+\right) \\
 s_3^*(2i+8) & = \left(1-\eta -\frac{D(1-D)}{1+D^2}+w^+\right)\left(1-2\beta +\beta D -\frac{\beta D^2 (1-D)}{1+D^2}-\beta D^3 w^-\right). \endalign $$
Now $w=\min\{w^+,w^-\}\geq 0$ else the minimum of these is certainly less than
$$\left(1-\eta -\frac{D(1-D)}{1+D^2}+w\right)\left(1-2\beta +\beta D -\frac{\beta D^2 (1-D)}{1+D^2}-\beta D^3 w\right)<\delta_{\infty}. $$
Hence we can assume that the sequence for $\gamma$ consists of blocks $A_1A_1'(A_1A_1A_1'A_1')^{l_i}$. Now if our block occurs inside a block
$$ ...0,1,0,-1 (0,1,0,1,0,-1,0,-1)^l0,1,0,-1,0(1,0,1,0,-1,0,-1,0)^m1,0,-1,0...$$ 
with $l>m$ then $w^-\leq 2D^{1+4l}\leq 2D^{4m+5}$ and
$-w^+\leq -2D^{4m+1}(1-D-D^2)$ and $s_1^*(2i+1)\leq \delta_{\infty}-2\beta D^{4m+4} ((1-\eta -D)(1-D-D^2)-\eta ) <\delta_{\infty}.$
Hence we need only worry about the equal lengthed blocks of $S_{k}$.
So the only values above $\delta_{\infty}$ come from these $\delta_{k}$, $k\geq 0$,
and  $\delta_{-1}$, plus  $\delta_{-2}$ when $b<2a$ and $\delta_{-2}>\delta_{\infty}$ (one tediously checks that this and $\delta_{-2}>\delta_{0}$ occur 
as stated).

Suppose now that $b$ is even. From the above 
we  need only consider blocks $A_0,A_2,C_0,C_2,C_4$ or their negatives.

We also note that if $t_{2i+1}=-2$ and $t_{2i}=0$ (or $t_{2i+2}=0$)
then
$$ \align s_1^*(2i)  & \leq \left(1-\beta +\frac{2D}{1-D}\right)\left(1-3\eta +\frac{4D-2\eta D}{1-D}\right)\\
  & \leq 1-3\eta-\beta+9D=u-2D(b-a-3-\eta+\beta),\endalign $$
is  less than $u$ for $b\geq a+4$ and 
$$ s_1^*(2i)\leq \left(1-\beta +\frac{2D}{1-D}\right)\left(1-3\eta +\frac{2D}{1-D}\right)\leq 1-3\eta-\beta+7D-6\eta D, $$
for $b=a+2$ enabling us to rule out blocks $C_0$.

If $t_{2i-1}=-2$ and $t_{2i}\leq 2$ then
$$ \align s_1^*(2i-1)  & \leq \left(1+\beta +\frac{2\beta D}{1-D}\right)\left(1-3\eta +\frac{4D-2\eta D}{1-D}\right)  \\
  & \leq 1-3\eta+\beta +D+6\beta D-2\eta D=u-2D(b-2a+1-\beta). \endalign $$
Hence for $b\geq 2a$ we can assume that either the expansion of
$\gamma$ eventually
has period $C_4$ 
or that  $t_{2i-1}=0$ for all $i$.

Suppose first that $b\geq 2a$.
Now infinitely many blocks of the form ${A'}_2{A'}_2$
or ${A'}_2A_0{A'}_2$
or ${A'}_2A_0A_0$ (or their negatives or reverses)
would lead to
$$ M^*(\alpha,\gamma)\leq \left(1-3\beta +\frac{2\beta D}{1-D}\right)\left(1-\eta +\frac{2 D^3}{1-D}\right)\leq \delta_{\infty,1}. $$
Hence if $\gamma$ does not have period $C_4$ or $A_0$ or $A_2{A'}_2$
then we can assume that its expansion consists of infinitely many $A_0$
separated by alternating $A_2$ and ${A'}_{2}$.
We suppose that $t_{2i}=0$, $t_{2i-2}=-2$, $t_{2i+2}=2$ and write
$d_{2i+1}^+=\frac{2\beta}{1+D}+w^+,$ $d_{2i-2}^-=-\frac{2\beta}{1+D}-w^-,$
so that
$$\align
s_3^*(2i+1)&=\left(1-\beta-\frac{2\beta}{1+D} -w^+\right)\left(1-\eta +\frac{2D^2}{1+D}+\eta Dw^-\right)\\
s_1^*(2i-2)&=\left(1-\beta-\frac{2\beta}{1+D} -w^-\right)\left(1-\eta +\frac{2D^2}{1+D}+\eta Dw^+\right)
\endalign$$
We can clearly assume that $w=\max\{w^+, w^-\}<0$ else the minimum of these will be at most
$$\left(1-\beta -\frac{2\beta}{1+D}-w\right)\left(1-\eta +\frac{2D^2}{1+D}
+\eta D w\right) 
   \leq \delta_{1,\infty} -w(1-\eta-\eta D)<\delta_{\infty,1}. $$ 
Hence the $t_{2i}=0$ occur in the 
middle of  blocks
$$ A_0 (A_2{A'}_{2})^mA_0(A_2{A'}_{2})^lA_0. $$
Now if $l>m$ then $w^-\leq -2\beta (1-2D)D^{2m}$ and $w^+\geq -2\beta D^{2l}\geq -2\beta D^{2m+2}$ giving
$$ s_{3}^*(2i+1)\leq \delta_{\infty,1} -2D^{2m+2}((1-2D)(1-3\beta)-\beta)<\delta_{\infty,1}.$$
Hence this leaves only the periodic case $A_0(A_2{A'}_2)^l$
and the result follows (after checking that $\delta_{0,4}>\delta_{1,\infty}$ iff $2a\leq b\leq 3a-6$).

Suppose next that $b=2a-2$, $a\geq 8$ and note the bound
$$\delta_{\infty,2}\geq 1-3\eta+\beta+D-2\eta D+2\beta D.$$
  If $t_{2i-1}=-2$, $t_{2i}=2$ and either $t_{2i+1}=-2$ or $t_{2i+1},t_{2i+2},t_{2i+3}=0,-2,0$
then
$$ s_1^*(2i-1)\leq \left(1-3\eta +\frac{4D-2\eta D}{1-D}\right)
\left(1+\beta -2\beta D+\frac{2\beta D^2}{1-D}\right) <\delta_{\infty,2},$$
and if we have blocks $t_{2i-2},t_{2i-1},t_{2i}=2,-2,2$ then
$$ s_1^*(2i-1)\leq \left(1-3\eta +\frac{2D}{1-D}\right)\left(1+\beta +\frac{2\beta D}{1-D}\right)\leq 1-3\eta +\beta -D+4\beta D <\delta_{\infty,2}, $$
and if we have blocks $t_{2i-2},...,t_{2i+1}=4,-2,2,0,$ with $t_{2i+2}\geq 0$ then
$$ \align s_3^*(2i) & \leq \left(1-3\beta +2D-4\beta D+\frac{2D^2}{1-D}\right)\left(1-\eta +\frac{2\eta D}{1-D}\right) \\
  & \leq 1-3\beta-\eta +5D-4\beta D+4\eta D^2 <\delta_{\infty,2}.\endalign  $$
Hence if we have infinitely many $t_{2i-1}=-2$ then either $\gamma$ has period $C_4$ or we have infinitely many blocks $t_{2i-2},...,t_{2i+4}=4,-2,2,0,-2,2,-4$.
Hence we have 
$$\align s_1^*(2i-1) & = (1-3\eta +\eta w^-)(1+\beta -2\beta D +2D^2-D^2w^+) \\
 s_3^*(2i+2) & = (1-3\eta +\eta w^+)(1+\beta -2\beta D +2D^2-D^2w^-) \endalign $$
where $w^+,w^-\leq (4\beta -2D)/(1-D)$ and $\min \{
s_1^*(2i-1),s_3^*(2i+2)\}\leq \delta_{\infty,2}$. Hence it remains only to 
deal with cases where $t_{2i-1}=0$ for all $i$ as for $b\geq 2a$.
Since $\delta_{0,5}\leq 1-\eta -3\beta +5D-4\beta D<\delta_{\infty,2}$ 
for $b=2a-2$, these only add $\delta_{0,1}$.

Suppose next that $b=2a-4$, $a\geq10$.
We note the upper bound
$$\delta_{\infty,3}\geq (1-\eta-2\eta D+2D)(1-3\beta+2D-4\beta D)\geq 1-3\beta-\eta +7D-4\eta D-10\beta D, $$
and rule out the configuration
$t_{2i-2},t_{2i-1}=-2,0$ 
with $t_{2i}\leq 0 $, since this leads to
$$ s_1^*(2i-2)\leq (1-3\beta +2D)\left(1-\eta +\frac{2D^2}{1-D}\right)\leq 1-\eta -3\beta+5D <\delta_{\infty,3}, $$
and $t_{2i-2},...,t_{2i+1}=2,-2,2,-2$,
since for $a\geq 10$ this leads to
$$\align  s_1^*(2i-1) & \leq (1-3\eta +2D)\left(1+\beta-2D+\frac{4\beta D-2D^2}{1-D}\right)\\
  & \leq 1+\beta-3\eta -3D +6\eta D+6\beta D<\delta_{\infty,3}. \endalign $$

Suppose that we have infinitely many blocks $t_{2i-2},t_{2i-1},
t_{2i}=-2,0,2$. 
If  $t_{2i-3}=0$, or $t_{2i-3},t_{2i-4}=2,-4$ 
then (4.1) with $d\geq 2\beta -
\frac{(2D-4\beta D+2D^2-2\beta D)}{1-D^2}$ 
(since we have excluded blocks $C_2C_2$ or $C_0$)
gives  $s_1^*(2i-2)$ or $s_3^*(2i-1)$ less than $\delta_{\infty,3}$. So we need only include
periods $A_0$ or $A_2C_2A_2'C_2'$ or $\gamma$ consisting solely
of blocks $C_2$ and $C_4$ with no blocks $C_2C_2$.
Hence either $\gamma$ has period $C_4$ or $C_2C_4$ or $\gamma$ consists
of blocks $...,b_{2i}|b_{2i+1},...$ of the form
$...,C_4(C_4C_2)^k|C_4C_4(C_2C_4)^lC_4,... $. Hence,
writing $d_{2i+1}^-=\frac{-2\eta}{1-D}+\frac{2D+4D^2}{1-D^2}+w^-$ and
 $d_{2i+4}^+=\frac{-2\eta}{1-D}+\frac{2D+4D^2}{1-D^2}+w^+$,
if $0\leq l<k$ we have $w^-\leq \frac{2D^{2k+1}}{1-D^2}\leq \frac{2D^{2l+3}}{1-D^2}$ and $w^+\geq \frac{2D^{2l+1}}{1-D^2}(1-D)$ and $s_2^*(2i+1)$ is
$$\align 
  & \left(1-\eta -\frac{2\eta D}{1-D}+\frac{2D+4D^2}{1-D^2}
+w^-\right)\left(1-3\beta +\frac{2D}{1-D}-\frac{(4\beta D+2\beta D^2)}{1-D^2}
-\beta Dw^+\right)\\
 & \leq \delta_{\infty,3}-\frac{2\beta D^{2l+2}}{1-D^2}((1-\eta)(1-D)-\eta)<\delta_{\infty,3}, \endalign$$
leaving only the periodic elements of $S_{k,3}$.

Suppose now that $b\leq 2a-6$ or $a=6,b=8$.
If $t_{2i+1}=0$ and $t_{2i}$ or $t_{2i+2}=\pm 2$ infinitely often then (4.1) with
$d\geq (2\beta -2D)/(1-D)$ gives $M^*(\alpha,\gamma)\leq \delta_{\infty,4}$.



Hence when $b\leq a+4\leq 2a-6$, or $a=6,b=8$, the only values above
$\delta_{\infty,4}$ come from periods $A_0$ or $C_2$. When $a+6\leq b\leq 2a-6$
we need also to consider $\gamma$ consisting of $C_2$ and $C_4$ with infinitely
many $C_4$. Suppose we have blocks $...,b_{2i-2} | b_{2i-1},...$
of the form $...C_4C_2^l|C_4C_2^mC_4...$ with $l>m\geq 0$ then, 
writing $d_{2i-1}^-=\frac{-2\eta+2D}{1-D}+w^-$ and $d_{2i-1}^+=4\beta -\frac{2D(1-\beta)}{1-D}+\beta w^+$, the bounds $w^-\leq 2D^{l+1}/(1-D)\leq 2D^{m+2}/(1-D)$
and $w^+\geq 2D^{m+1}$ give
$$
\align
s_2^*(2i-1) & = \left(1-\eta +\frac{2D(1-\eta)}{1-D}+w^-\right)\left(1-3\beta+\frac{2D(1-\beta)}{1-D}-\beta w^+\right) \\
 & \leq \delta_{\infty,4}-2\beta D^{m+1}((1-\eta)-\eta)\leq \delta_{\infty,4}.\endalign $$
Hence it is enough to consider the periodic elements in $\delta_{k,4}$, $k\geq 0$.   

It remains to deal with the left-over cases $(a,b)=(4,6)$, $(6,10)$ and $(8,12)$.

When $(a,b)=(4,6)$ we can rule out blocks 
$t_{2i-3},...,t_{2i}=0,-2,0,-2$ else
$$s_1^*(2i-2)\leq \left(1-3\beta +\frac{2\beta D}{1-D}\right)\left(1-\eta -2D+\frac{2\eta D-2D^2}{1-D}\right)<\delta_{\infty,1}$$
and $t_{2i-2},...,t_{2i+3}=-2,0,0,0,-2,0$, else
$$s_1^*(2i+1)\leq \left(1-3\beta +\frac{2\beta D}{1-D}\right)\left(1-\eta -2D^2+2\eta D^2\right)<\delta_{\infty,1}.$$
Hence if we have blocks $t_{2i-1},t_{2i}=-2,2,$
with $t_{2i+1}=-2$, or with
$t_{2i+1}=0$ and  $t_{2i+2}\leq 0$,
then
$$ s_1^*(2i-1)\leq \left(1-3\eta+\frac{2D}{1-D}-\frac{2\eta D^2}{1-D^2}\right)\left(1+\beta +\frac{2\beta D^2}{1-D}-\frac{2D^3}{1-D^2}\right) <\delta_{\infty,1}. $$
So  either $\gamma$ has period $A_2C_2$ or we can assume that $t_{2i-1}=0$ for all $i$ (in which case the proof follows exactly as for $b\geq 2a$). 

For $(a,b)=(6,10)$ observe that we can have no blocks $t_{2i},t_{2i+1},t_{2i+2}=-4,2,-4$ else
$$ s_1^*(2i-1)\leq \left(1-5\beta+\frac{2D-2\beta D}{1-D}\right)\left(1+\eta -4D+\frac{2\eta D-2D^2}{1-D}\right)<\delta_{\infty,7}, $$
and if we have $t_{2i-1}=-2$ then $t_{2i-2}$ or $t_{2i}=4$ else
$$ s_1^*(2i-1)\leq \left(1-3\eta +\frac{2D}{1-D}\right)\left(1+\beta +\frac{2\beta D}{1-D}\right) <\delta_{\infty,7}. $$
Likewise we can rule out blocks 
$t_{2i-1},...,t_{2i+1}=0,-2,0$ and $t_{2i+2}\leq 0$
else
$$ s_1^*(2i)\leq \left(1-3\beta +\frac{2\beta D}{1-D}\right)\left(1-\eta +\frac{2D^2}{1-D}\right) <\delta_{\infty,7}.$$
Hence if $\gamma$ has infinitely many $t_{2i-1}=\pm 2$ but not period
$C_2C_4$ then we must have infinitely many blocks $t_{2i-2},...,t_{2i+1}=-4,2,-2,0$ 
and
$s_1^*(2i-2)\leq \delta_{\infty,7}$. If $t_{2i-1}=0$ for almost all
$i$ then all but period $A_2A_2'$ and $A_0$ have been ruled out.

For $(a,b)=(8,12)$ we first rule out
$t_{2i-1},...,t_{2i+3}=-2,2,-2,2,-2$ else
$$ s_1^*(2i+1)\leq \left(1-3\eta +2D-2\eta D+\frac{2D^2}{1-D}\right)\left(1+\beta -2D +\frac{2\beta D}{1-D}\right) < \delta_{\infty,6}. $$
If $t_{2i-1},t_{2i}=0,-2$
then $t_{2i+1},t_{2i+2}=2,-2$ 
else
$$ s_1^*(2i-1)\leq \left(1-\eta +\frac{2D}{1-D}\right)\left(1-3\beta +\frac{2\beta D}{1-D}\right) <\delta_{\infty,6} $$
and $t_{2i-2},t_{2i-3}=2,-2$ else
$$ s_1^*(2i-1)\leq \left(1-\eta +\frac{2D^2}{1-D}\right)\left(1-3\beta +\frac{2D-2\beta D}{1-D}\right) <\delta_{\infty,6}. $$
Likewise if $t_{2i-1},t_{2i-2}=0,-2$ then $t_{2i-4},...,t_{2i+1}=-2,2,-2,0,2,-2$. Hence blocks $A_2'$ occur inside blocks $C_2A_2'C_2'$ and either 
$\gamma$ has period  $A_2C_2A_2'C_2'$ or we have infinitely many blocks
$b_{2i-1},...,b_{2i+4}=C_2'A_2C_2$ with $b_{2i+5}b_{2i+6}=C_2$ and,
writing
$$ d_{2i+3}^+,-d_{2i-2}^-=\frac{2\beta +2\beta D-2D}{1+D^2}+w^+, w^-$$ 
we have
$$ \align s_1^*(2i+3)=  & \left(1+\beta -2D\frac{(1-\beta +\beta D)}{1+D^2}+w^+\right)  \\
  &  \;\; \left(1-3\eta +2D-2D^2+2\eta D^2-2D^3\frac{(1-\eta+D)}{1+D^2} -\eta D^2 w^-\right), \endalign  $$
with a similar expression for $s_{3}^*(2i-2)$ with $w^+$ and $w^-$ interchanged.
Clearly then $w=\min\{w^+,w^-\}>0$ else the minimum of $s_1^*(2i+3)$ and $s_3^*(2i-2)$ will be less than 
$\delta_{\infty,6}+w((1-3\eta)-\eta D^2(1+\beta))<\delta_{\infty,6}.$ 
Now if this block $C_2'A_2C_2$ occurs inside a block
$$ ...C_2C_2 (A_2'C_2'A_2C_2)^l A_2'C_2'C_2'A_2C_2C_2 (A_2'C_2'A_2C_2)^k A_2'C_2'C_2'... $$
with $k>l\geq 0$ or
$$ ...A_2'C_2' A_2C_2C_2 (A_2'C_2'A_2C_2)^kA_2'C_2'C_2'... $$
with $k\geq 0$, then $w^+\leq 2D^{4+4k}$ and $w^-\geq 2D^{4+4l}(1-2\beta)$ or 
$w^-\geq 2D(1-2\beta)$
respectively, and $w^-\geq 2D^{1+4k}(1-2\beta)$ gives
$$ s_1^*(k)\leq \delta_{\infty,6}-2\eta D^{3+4k}(1-2\beta)+2D^{4+4k}\leq \delta_{\infty,6}, $$
leaving only the periodic elements of $\delta_{k,6}$. $\blacksquare$


\head 5. Proof of Theorem 2 (odd $a\geq 3$) \endhead

We first exclude  the cases $a=3,b=4,5$ or 6.
In all the remaining cases  the $\gamma$ of interest have
 $t_{i}\neq a_i$ for almost all $i$. Now if $t_i=a_i$ infinitely
often then from (1.2) we have $M^*(\alpha,\gamma)\leq \eta$.
Writing
$$\align \delta_{\infty,1} & = \left(1-2\eta+\eta v+\frac{2D}{1-D}\right)\left(1-\beta +v+\frac{2\beta D}{1-D}\right) = \eta + E_1 \endalign $$
with
$$ E_1=\eta (a-3)(1-\beta)+ \beta D-\frac{2D^2}{1-D} +v\left(1-\eta+D+\frac{4D^2}{1-D}\right)+\eta v^2 +\frac{4\beta D^2}{(1-D)^2}. $$
When $m\geq 1$ we plainly have $v>0$ and $E_1>0$, when $m=0$ and $a\geq 5$ we have $E_1\geq \eta 2(1-\beta)-D >0$.
Hence $\delta_{\infty,1}\geq \eta$ and we can rule out $b_i=a_i-1$ when
$r\geq a+3$ and $b\geq 7$.
For $r\leq a+1$ we use that
$$ \align \delta_{\infty,3} & =\left(1-\eta v-\frac{2D}{1-D^2}\right)\left(1-3\beta -v -\frac{2\beta D^2}{1-D^2}\right) \\
  & > (1-Dm-2D)(1-3\beta -m\beta +D-m\beta D) \\
  & > \eta + E_2  \tag 5.1 \endalign $$
with
$$ E_2=1-\eta -3\beta -m\beta -D-Dm=\beta \left( m(a-2-\eta)+r\left(1-\frac{1}{a}\right) -3-2\eta \right). $$
If $m=0$ and $r=a+1$ then $E_2=\beta (a-1/a-3-2\eta)>0$ for $a\geq 5$. If $m\geq 1$
and $a\geq 5$ then $E_2\geq \beta (a-3-2/a-3\eta)>\beta >0$. 
For $a=3$ ($r=2$ or 4) we have $E_2>\beta (\frac{1}{2}m-\frac{8}{3})>0$ for
$m\geq 6$. One checks numerically that $\delta_{\infty,3}>\eta$ also holds
for $b=3m+4, m=1,...,5$ and $b=3m+2$, $m=3,4,5$, with $\delta_{\infty,7}>\eta$
for $a=3,b=8$. Hence we can assume that $t_i\neq a_i$ for almost all $i$.

We next show that (other than when $r=2,4$, $m=1$) the crucial $\gamma$ have
$t_{2i+1}=\pm 1$. From (1.3) we know that if $|t_{2i+1}|\geq 3$ then $M^*(\alpha,\gamma)<1-3\eta +D$. 
The rough lower bound 
$$ \delta_{\infty,1} >1-2\eta -\beta +3D $$
is enough to rule out this possibility when $r\geq a+3$. From
(5.1) we similarly obtain
$$ (1-3\eta +D) \leq \delta_{\infty,3} +\beta E_3$$
with
$$ E_3 =3+2\eta -m (2-\eta) -\frac{3r}{a}. $$
Hence if $r=a+1$ we have $E_3<(2\eta -3/a)<0$ and if $m\geq 2$ 
we have $E_3< 4\eta -1 - \frac{3r}{a} <0$. 

Hence we only need consider $t_{2i+1}\neq \pm 1$ when $m=1$
and  $r\leq a-1$.  Now for  $m=1$ if $|t_{i}|\geq 5$
infinitely often  then $M^*(\alpha,\gamma)\leq 1-5\beta +D<\delta_{\infty,3}-\beta +3D <\delta_{\infty,3}$. Hence for $m=1$, $r\leq a-1$ we can assume
that $t_i=\pm 1,\pm 3$.
When $m=1$, $r\leq a-1$ we note the rough bound
$$ \align \delta_{\infty,5} & >(1-D+\eta D-D^2)(1-4\beta +D-\beta D)\\
 & > 1-4\beta +3\beta D. \endalign  $$
Now if for example  $t_{2i+1}=-3$ and $t_{2i+2}\leq 1$ (likewise if $t_{2i}\leq 1$) then
$$\align s_1^*(2i+1) & \leq \left(1-4\eta +\frac{3D+3\eta D}{1-D}\right)\left( 1-\beta +\beta +\frac{3D+3\beta D}{1-D}\right) \\
 & \leq 1-4\eta +\frac{3D+3\eta D}{1-D} +(1-3\eta)\frac{(3D+3\beta D)}{1-D} \\
 & < 1-4\eta +6D =1-4\beta +D\left(6-\frac{4r}{1+D}\right) <\delta_{\infty,5}. \endalign $$
Hence if $t_{2i+1}=3\delta$, $\delta=\pm 1$ we can assume that
$t_{2i},t_{2i+2}=-3\delta. $
Similarly, if $t_{2i+2}=-3$ and $t_{2i+1}\leq -1$ (likewise for $t_{2i+3}$) then
$$ \align s_{1}^*(2i+1) & \leq (1-\eta -\eta +3D+3\eta D)(1-\beta -3\beta +3D) < 1-4\beta <\delta_{\infty,5}, \endalign $$
and we assume that  $t_{2i+2}=-3\delta$ implies 
$t_{2i+1},t_{2i+3}=\delta$ or $3\delta$.
Now if $t_{2i+1},t_{2i+2},t_{2i+3}=-3,3,-1$
then
$$ \align 
s_2^*(2i+1) & \leq \left(1+\eta -3\eta +\frac{3D-\eta D}{1-D}\right)\left(1+\beta -3\beta +\frac{D(1+\beta)}{1-D}\right) \\
& \leq (1-2 \eta +3D)(1-2\beta +D +2\beta D) \\
 & < 1-2\beta -2\eta +8D-2\eta D-4\beta D \\
 & < \delta_{\infty,5} -2\beta \frac{(r-4)}{a}.\endalign $$
Hence when $r\geq 4$ we can assume that either $t_{2i+1}=\pm 1$
for almost all $i$ or that $\pm \gamma$ is in $S_{-5}$.
Moreover since
$$ \align \delta_{-5} & \leq (1-2\eta +3D)(1-2\beta+3D) \\
   & < 1-2\beta -2\eta +10D < \delta_{\infty,5}-2\beta \left(\frac{r}{a}-5\eta\right), \endalign $$
we do not even  need to include $\delta_{-5}$ for $r\geq 6$.
For $r=4$ and  $5\leq a\leq 11$ one checks numerically that $\delta_{\infty,5}>\delta_{-5}$.
For $r=4$ and $a\geq 13$ we have
$$ \delta_{-5}-\delta_{\infty,5}=\frac{2D}{1-D}\left(1-11\beta -\frac{6D(1-D)}{1+D}+\frac{5D (1-\eta)(1-\beta)}{1-D}\right)>0,$$
with $\delta_0>1-2\eta >\delta_{-5}$.

Suppose now that $b=a+2$. If $a=5$ or  $7$ then we can still assume
that $t_{2i+1}=\pm 1$ since if $t_{2i+1}=-3$
then
$$ s_{1}^*(2i+1)\leq \left(1-\eta-3\eta+\frac{3D-\eta D}{1-D}\right)\left(1-\beta+\frac{3\beta-D}{1-D}\right)<\delta_{\infty,8}. $$ Hence we assume
that $a\geq 9$. We note the lower bound
$$\align
\delta_{\infty,9} & \geq \left(1-2\eta+3D-3\eta D\right)\left(1-2\beta +D-\beta D\right) \\
 & > 1-2\beta-2\eta+8D-5\eta D-7\beta D+10D^2 >1-4\beta +4D-12\eta D. \endalign $$
Now if $t_{2i+2}=-3\delta$, $\delta=\pm 1$, then we can assume
that at least one of $t_{2i+1},$ $t_{2i+3}=3\delta$, since
if
$t_{2i+2}=-3,$ $t_{2i+1}=t_{2i+3}=1$ then
$$ \align 
s_1^* (2i+1) & \leq \left(1-\eta +\frac{\eta +D}{1-D}\right)\left(1-\beta -3\beta+\frac{D+\beta D}{1-D}\right) \\
  & < 1-4\beta +\frac{D+\beta D}{1-D}+\frac{D+\eta D}{1-D}(1-3\beta)
\\
 & <1-4\beta +2D <\delta_{\infty,9}. \endalign $$
Similarly if $t_{2i+2}=-3\delta$ and $t_{2i+1}=\delta$ then $t_{2i}=\delta$, since if
$t_{2i+2},t_{2i+1}=-3,1$, $t_{2i}\leq -1$ then
$$\align 
s_1^* (2i+1) & \leq \left(1-\eta+\eta -D+\frac{\eta D+D^2}{1-D}\right)
\left(1-\beta-3\beta +\frac{3D-3\beta D}{1-D}\right)\\
 & <(1-D+2\eta D)(1-4\beta+3D-2\beta D)  \\
  & <1-4\beta+2D+2\beta D+2\eta D <\delta_{\infty,9}. \endalign $$
(Likewise $t_{2i+4}=\delta$ if $t_{2i+2}=-3\delta$ and $t_{2i+3}=\delta$). 
Also if $t_{2i+1}=\delta$ and $t_{2i+1+\delta'}=-\delta$ with $\delta$ and $\delta'=\pm1 $ then $t_{2i+1-\delta'}=-\delta$ or $-3\delta$, since if $t_{2i+1},t_{2i+2}=
-1,1,$ $t_{2i}\leq -1$ then
$$\align 
s_{1}^*(2i) & \leq \left(1-\beta-\beta +\frac{D+\beta D}{1-D}\right)\left(1-\eta -\eta +\frac{D+\eta D}{1-D}\right) \\
  & < (1-2\beta)(1-2\eta)+\frac{D}{1-D}((1+\beta)(1-\eta)+(1+\eta)(1-2\beta)) \\
 & < 1-2\beta -2\eta +6D <\delta_{\infty,9}. \endalign $$
So if $t_{2i+1}=-1$, $t_{2i+2}$ (or $t_{2i}$) 
$=-1$ infinitely often then  $s_{1}^*(2i+1)$ is at most
$$ \left(1-\eta-\eta+3D-3\eta D+3D^2-\eta D^2+\frac{-D^3+\eta D^3}{1-D}\right)\left(1-\beta -\beta +\frac{D-D\beta}{1-D}\right). $$
Thus  if $\pm \gamma$ has $t_{2i}=-3$ infinitely 
often and $\pm \gamma$ is not in $S_{-5}$ we must have infinitely many
blocks $t_{2i+1},t_{2i+2}$ or $t_{2i-1},t_{2i-2}$ of the form $1,1$ and hence $M^*(\alpha,\gamma)\leq \delta_{\infty,9}$. This only leaves the $\gamma$ with $t_{i}=\pm 1$, for almost all $i$ with $t_{2i+1},t_{2i+1\pm 1}$
 of the form $-1,-1$ or $1,1$ at most finitely often; which leaves just the $\pm \gamma$ in $S_0$.
For $a\geq 11$ we have
$$\delta_{-5}>\delta_{\infty,9}-2\eta D^2+2D(1-\beta)(1-2\eta)>\delta_{\infty,9}.$$
One checks numerically that this is also true for $a=9$.

Thus in the remaining cases we can assume that $t_{2i+1}=\pm 1$ for all $i$.
Suppose that $\gamma$ and $1-\alpha-\gamma$ satisfy
$$\align t_{2i+1}=-1 & \Rightarrow (m-S) \leq t_{2i},t_{2i+2}\leq (m+2+R) 
\tag 5.2 \\ 
 t_{2i+1}=1 & 
\Rightarrow -(m+2+R)\leq t_{2i},t_{2i+2}
\leq -(m-S) \endalign $$
for almost all $i$, with $R$ and $S$ minimal.
 
Suppose first that $R\geq S$ and that $\gamma$ (or $1-\alpha -\gamma$) has
$ t_{2i+1}=1$, $t_{2i+2}=-(m+2+R)$
infinitely often.
Hence
$$d_{2i+1}^-\leq \eta + \frac{ (-m+S)D +\eta D}{1-D}
 $$
$$ d_{2i+1}^+\leq -(m+2+R)\beta +\frac{D}{1-D} + \frac{(-m+S)\beta D}{1-D}, $$
and
$$\align  s_1^*(2i+1) & \leq \left( 1-\eta v + \frac{RD}{1-D}\right) \left(1
-3 \beta -v  -\frac{\beta R(1-2D)}{1-D}\right) \\
 & = \delta_{\infty,5} -\frac{R\beta}{1-D}\left(1-\eta+D+2\eta D v +\frac{RD (1-2D)}{1-D}\right) \\
 & < \delta_{\infty,5} -\frac{2\beta}{1-D}(1-\eta +D), \endalign $$
for $R\geq 2$. Since
$$ \align \delta_{\infty,5} & = \delta_{\infty,3} +\frac{2\beta /a}{1-D}\left( 1-3\beta+D+\beta D-\frac{m\beta (1+D^2)}{1-D}+\frac{D^2(1-2/b+D)}{1-D}\right) \\
  & < \delta_{\infty,3} +\frac{2\beta/a}{1-D}, \endalign $$
and 
$$ \align 
\delta_{\infty,5} & = \delta_{\infty,1} -\frac{2\beta}{1-D}\left(1-\frac{r}{a}(1-D)+\eta -D-\frac{2\eta D}{1-D}+\frac{2(m+1)D^2}{1-D}\right) \\
 & < \delta_{\infty,1} -\frac{2\beta}{1-D}(-1+\eta +D -2\eta D), \endalign $$
we immediately obtain $M^*(\alpha,\gamma)<\delta_{\infty,3}$ and $\delta_{\infty,1}$ for $R\geq 2$.
Hence we can assume that $S\leq R\leq 0$ and (since plainly $R+S
\geq -2$)
the only possibilities are $R=S=0$ or $R=0,S=-2$.
 
Similarly suppose that $R<S$ and $t_{2i+1}=-1$, $t_{2i+2}=(m-S)$.
Then
$$ d_{2i+1}^-\leq -\eta +(m+2+R)D 
+\frac{\eta D}{1-D}+\frac{(-m+S)D^2}{1-D},$$
$$d_{2i+1}^+\leq (m-S)\beta +
\frac{D}{1-D} +\frac{(-m+S)\beta D}{1-D}, \tag 5.3  $$
bounding $s_1^*(2i+1)$ by
$$ \align & \left(1-2\eta +\eta v+\frac{2D}{1-D}+\frac{D}{1-D}(S-2+2\eta-2mD)\right) \\
  & \hskip 5ex \left(1-\beta+v+\frac{2\beta D}{1-D}-\frac{\beta}{1-D}(S(1-2D)+2D+2mD-2\eta)\right) \\
 &= \delta_{\infty,1} -\frac{2D}{1-D}\left( \frac{1}{2}S(a-3)+\frac{1}{2}(S-2)U_1 + U_2\right), \tag 5.4  \endalign  $$
where (since $m\beta <\eta$) 
$$\align U_1 & = S\frac{\beta (1-2D)}{1-D}+\frac{2D}{1-D}\left(2-\beta -3D-m\beta (1-2D)\right)>0, \\
  U_2 & = \frac{2D(1-3D)}{1-D} +\eta +\frac{4D^2(m\beta)}{1-D} +(m\beta)(2-3\eta -D)+2(m\beta)^2 \eta>0. \endalign $$ 
Hence $M^*(\alpha,\gamma)<\delta_{\infty,1}$ for
$S\geq 2$. For $r=a+1$ we have 
$$ \delta_{\infty,1}-\delta_{\infty,3} =\frac{2D}{1-D}W, $$
with
$$\align  W &
=1-\eta -2\beta +D +\frac{2\eta D}{1-D^2}+\frac{4\beta D}{1+D}-\frac{2D^2(3+D-\beta D-D^2)}{(1+D)(1-D^2)}  \\
  & < 1-\eta +2\eta D-\frac{4}{3}\beta .\endalign $$
When $a\geq 5$ the bound $\frac{1}{2}S(a-3)\geq 2$ gives $M^*(\alpha,\gamma)<\delta_{\infty,3}$ for $S\geq 2$. When $a=3$ we have 
$$\align
U_2 &  =1-\eta +D -2\eta D +\frac{2D(1-3D)}{1-D}-\frac{r}{3}\beta (1-2D)+2\eta (m\beta)^2 +\frac{4D^2 (m\beta)}{1-D} \\
 & > 1-\eta +3D-2\eta D-\frac{4}{3}\beta > W \endalign $$
and still obtain $M^*(\alpha,\gamma)<\delta_{\infty,3}$. 
Writing
$$ \delta_{\infty,1}=\delta_{\infty,5} +\frac{2D}{1-D}W_1, $$
where
$$ W_1=(a-r+1)-2\beta \left(1-\frac{r}{a}\right)-\frac{2D}{1-D}\left(1-\eta -\beta \left(1-\frac{r}{a}\right)\right), $$
it is clear that for $4\leq r\leq a-1$ the rough bounds $W_1<(a-3)$ and 
$\frac{1}{2}S (a-3)\geq (a-3)$ are enough to give $M^*(\alpha,\gamma)<\delta_{\infty,5}$ when $S\geq 2$. 
Hence when $r\geq 4$ we must have $S\leq 0$, so we can assume that
$$ t_{2i+1}=\delta \Rightarrow t_{2i+1\pm 1}=-m\delta, \hbox{ or} -(m+2)\delta, $$
for $\delta=1$ or $-1$, so that the coefficients in the expansion of $\gamma$ consists of blocks $B_m$, $B_n$, $B_m'$ or $B_n'$.

When $r=2$ and $a\geq 5$ the bounds $W_1<(a-1)$ and $\frac{1}{2}S(a-3)\geq 2a-6$ 
give $M^*(\alpha,\gamma)<\delta_{\infty,5}$ for $S\geq 4$.  Thus we can assume that $S\leq 2$ and hence that  
$$ t_{2i+1}=\delta \Rightarrow t_{2i+1\pm 1}= -(m-2)\delta,-m\delta, \hbox{ or} -(m+2)\delta , \tag 5.5 $$
With a bit more work we show that this is also true for $a=3, b=3m+2, m\geq 2$:
For $a=3$, $r=2$ observe that $(m\beta)=\eta-\frac{2}{3}\beta$ and for $S\geq m$
$$ U_1\geq \eta +\frac{D}{1-D}
\left(1-\frac{5}{9}\beta +\frac{4}{3}D^2-\frac{41}{9}\beta D\right) >\eta $$
and
$$ \align U_2 & =1-\eta +D +\frac{2D}{1-D}((1-D)^2-\eta -\frac{4}{3}\beta D)-\frac{2}{3}\beta (1-2D) +2\eta (m\beta)^2 \\
  & > 1-\eta+D-\frac{2}{3}\beta > -4\eta +2-\frac{2}{3}\beta \endalign $$
while
$$ W_1=2-\frac{2}{3}\beta -\frac{2D}{1-D}(1-\eta-\frac{1}{3}\beta) <2-\frac{2}{3}\beta. $$
Hence $M^*(\alpha,\gamma)<\delta_{\infty,5}$ for $\frac{1}{2}(S-2)\geq 4$.
Thus it remains to consider $S<m$ or the odd cases $S=8, m\leq 8$, $S=6,m\leq 6$ or $S=4, m\leq 4$. Notice that $t_{2i}=(m+2+R)$ infinitely often implies that $M^*(\alpha,\gamma)\leq 1-(m+2+R)\beta +D$, so that we can 
assume that $(m+2+R)< (1+D-\delta_{\infty,5})/\beta $ is at most $12,11,10,7,6,5,4$ as $m=8,7,6,5,4,3,2$. Using this bound for $(m+2+R)$ in (5.3) gives
numerically $M^*(\alpha,\gamma)<\delta_{\infty,5}$ for the oddments.
Thus it remains only to deal with $4\leq S<m$, $R<S$. In particular (5.2)
implies that the $t_{i}$ are alternately greater or less than $0$, so that the $t_{2i+1}$ are constant, improving (5.4) to
$$\align 
M^*(\alpha,\gamma) & \leq \left(1-2\eta+\eta v+\frac{(R+2)D}{1-D}\right)\left(1-\beta +v-S\beta +\frac{(R+2)\beta D}{1-D}\right) \\
  & = \delta_{\infty,5} + DW_2 +D(R+2-S)v , \endalign $$
with
$$ \align
W_2 & = 2-\frac{2}{3}\beta -2v-S(1-\beta) +\frac{(R+2)}{1-D}\left(1-2D+2Dv -S\beta +\frac{D\beta (R+2)}{1-D} \right) \\
 & < (4+R-S)(1-\beta). \endalign $$
Hence $M^*(\alpha,\gamma)<\delta_{\infty,5} $ for $R\leq S-4$, so we can 
assume that $S\geq 4$ and
$ R= (S-2)\geq \frac{1}{2}S.$
But $t_{2i+1}=1$, $t_{2i+2}=-(m+2+R)$ then gives 
$$ \align s_1^*(2i+1) & \leq \left(1-\eta v +\frac{SD}{1-D}\right)\left(1-3\beta -v -\frac{1}{2}S\beta +\frac{S\beta D}{1-D}\right) \\
  & = \delta_{\infty,5} -\frac{1}{2}\frac{SD}{1-D}\left(\left(3-4\beta\right)\left(1-\eta v+\frac{SD}{1-D}\right)-2(1-3\beta -v)\right) \\
  & <\delta_{\infty,5}.  
\endalign $$
So we can assume that $S\leq 2$ and (5.5) holds.

Suppose now that $r\geq a+3$ with $m\geq 1$. From the previous discussion
we can assume that $t_{2i+1}=-1$ and $t_{2i+2}=m$ or $(m+2)$ for almost all $i$.

Now if $t_{2i+2}=m$ then
$$\align
s_1^*(2i+1) & \leq \left(1-\eta -\frac{\eta}{1-D}+\frac{(m+2)D}{1-D}\right)\left( 1-\beta +m\beta -\frac{D}{1-D} +\frac{(m+2)\beta D}{1-D}\right) \\
 & =\delta_{\infty,1}. \endalign $$
Hence either $\pm \gamma$ is in $S_{-2}$ or $M^*(\alpha,\gamma)\leq \delta_{\infty,1}$. For  $r\geq a+3$,
$$\delta_{-2}=\delta_{\infty,1}+2\beta \left( \frac{r-(a+2)}{a}+\eta v +\frac{2D}{1-D}-\frac{2D}{a}\right) > \delta_{\infty,1}.$$

Similarly suppose that $4\leq r \leq a-1$ and that $\gamma$ is not in $S_{-5}$ when $r=4$ and $a\geq 13$. Then from the above we can assume that $\gamma$ has
$t_{2i+1}=1$, $t_{2i+2}=-m$ or $-(m+2)$ for almost all $i$.  If $t_{2i+1}=-(m+2)$ then
$$\align 
s_{1}^*(2i+1) & \leq \left(1-\beta -(m+2)\beta +\frac{D-m\beta D}{1-D}\right)\left(1-\eta +\frac{\eta-mD}{1-D}\right)  =\delta_{\infty,5}. \endalign $$
Hence either $\pm \gamma$ is in $S_0$ or $M^*(\alpha,\gamma)\leq \delta_{\infty,5}$. Since $v\leq \eta$ and $r\leq a-1$,
$$ \delta_{0}=\delta_{\infty,5}+ 2\beta \left( \frac{a-r}{a}-\eta v\right)>\delta_{\infty,5}. $$

Suppose that $r=a+1$ with $m\geq 1$ then again we can assume that
$t_{2i+1}=1$, $t_{2i+2}=-m$ or $-(m+2)$ for almost all $i$.
If $t_{2i+2}=t_{2i}=-m$ then $s_3^*(2i+1)$ is at most
$$\align 
 &\left(1-\eta-\frac{\eta}{1-D} +mD+\frac{(m+2)D^2}{1-D}\right)\left(1-\beta +m\beta -\frac{D}{1-D}+\frac{(m+2)\beta D}{1-D}\right) \\
  &\;\;\;\;\; =\delta_{-2}-\frac{2\beta D(1+3D)}{1-D}+\frac{2\beta D}{a}-4Dv<\delta_{-2}<\delta_{\infty,3}. \endalign $$
Hence we can assume that the expansion of $\gamma$ eventually consists of $B_m$s
and $B_n$s with no consecutive $B_m$s. Thus either the $B_m$s and $B_n$s eventually
alternate and $\gamma$ is in $S_{-1}$ or $\gamma $ contains infinitely
many consecutive $B_n$s, $t_{2i+1\pm 1}=-(m+2)$ and
$$\align d_{2i+1}^+ & \leq  -(m+2)\beta +\frac{D}{1-D}-\frac{m\beta D}{1-D^2}-\frac{(m+2)\beta D^2}{1-D^2} \\ 
d_{2i+1}^- & \leq \frac{\eta}{1-D} -\frac{(m+2)D}{1-D^2}-\frac{mD^2}{1-D^2}
\endalign $$
giving  $s_1^*(2i+1) \leq \delta_{\infty,3}$.
Clearly $\delta_{-1}>\delta_{\infty,3}$ for $r=a+1$.

It remains to consider $m=0$ or $r=2$, $m\geq 2$.

Suppose first that $m=0$ and $r\geq a+3$. By the previous analysis we can assume that $t_{2i+1}=\delta_{i}$
and $t_{2i+1\pm 1}=0$ or $-2\delta_{i}$
with $\delta_{i}=\pm 1$. Hence if $\pm \gamma$ is not in $S_{-2}$ we must have
$t_{2i+2}=0$ infinitely often. Now if $t_{2i+2}=0$ we can assume that $t_{2i+1}\neq t_{2i+3}$, since
$t_{2i+1}=t_{2i+3}=-1$ would give
$$
s_1^* (2i+2) \leq \left(1-\eta -\eta +\frac{2D-\eta D}{1-D}\right)\left(1-\beta -D+\frac{2\beta D-D^2}{1-D}\right) < \delta_{\infty,2} . $$
We assume that $t_{2i+1},t_{2i+2},t_{2i+3}=1,0,-1$ giving
$$\align 
s_1^*(2i+2) & =(1-\eta -\eta +\eta d_{2i+3}^+)(1-\beta +D+Dd_{2i}^-) \\
s_3^*(2i+1) & = (1-\eta -\eta -\eta d_{2i}^-)(1-\beta +D-Dd_{2i+3}^+). \endalign $$
Hence, writing $d^*=\max\{d_{2i+3}^+,-d_{2i}^-\}$ we plainly have
$$ \min\{s_1^*(2i+2),s_3^*(2i+1)\} < (1-2\eta +\eta d^*)(1-\beta +D-Dd^*) \leq d_{\infty,2}^*, $$
since $d^*\leq (2\beta -D)/(1-D)$.
The claim follows on checking that
$$\delta_{-2}=\delta_{\infty,2} +\frac{2\beta}{(1-D)^2}\left(\frac{r-(a+3)}{a}+D\left( 6-\frac{3}{a}-2\eta+\eta D -D\right)\right) >\delta_{\infty,2}. $$

Suppose that $b=a+1$, $a\geq 5$. We note the rough bound
$$\delta_{\infty,4} >(1-\eta D)(1-3\beta +D) > 1-3\beta +D-\eta D. $$
From above we have
$t_{2i+1}=\delta_i, $ $t_{2i+1\pm 1}=0$ or
$-2\delta_i$, $\delta_{i}=\pm 1$ for almost all $i$. 
If $t_{2i+1\pm 1}=-2$ then
$$\align s_1^*(2i+1) & \leq \left(1-\eta +\eta -2D+\frac{\eta D}{1-D}\right)\left(1-\beta -2\beta +\frac{D}{1-D}\right)  \\
  & \leq 1-3\beta +D-\eta D-2D(1-\eta-3\beta)<\delta_{\infty,4}. \endalign $$
So we can assume that $\gamma$ does contain consecutive $t_{2i}=t_{2i+2}=-2\delta_i$. Similarly we rule out the possibility of $t_{2i}=t_{2i+2}=0$
infinitely often when $\gamma$ is not in $S_{-4}$, since otherwise 
we can assume that
$\pm \gamma$ has blocks $t_{2i+1}=-1$, $t_{2i}=t_{2i+2}=0$ and 
either $t_{2i+3}=-1$ and
$$\align
s_1^*(2i+1) & \leq \left(1-\eta -\eta +\frac{\eta D}{1-D}\right)\left(1-\beta -D +\frac{2\beta D-D^2}{1-D}\right) \\
  & < 1-2\eta -\beta +D+3\eta D+2\beta D =1-3\beta -D+5\eta D<\delta_{\infty,4}, \endalign $$
or $t_{2i+3}=1$ and $t_{2i+4}=-2$ and
$$\align
s_1^* (2i+3) & \leq \left(1-\eta +\eta -\eta D+\frac{\eta D^2}{1-D}\right)\left( 1-\beta -\frac{2\beta}{1-D^2} +\frac{D}{1-D}\right)<\delta_{\infty,4}.\endalign $$
Hence apart from $S_{-4}$ it is enough to consider $\gamma$ such that the 
$t_{2i}$ alternate between $0$ and $-2\delta_i$. 
Thus if $\pm \gamma$ is not in $S_{-1}$ we can assume that $\gamma$ contains
infinitely many blocks $t_{2i-1},...,t_{2i+5}=1,-2,1,0,-1,2,-1. $
Hence we have
$$\align
s_1^* (2i) & =(1-\beta -2\beta +D+Dd_{2i-2}^-)(1-\eta +\eta -\eta D+2D^2-\eta D^2 +\eta D^2 d_{2i+5}^+) \\
s_3^* (2i+3) & =(1-\beta -2\beta +D-Dd_{2i+5}^+)(1-\eta +\eta -\eta D+2D^2-\eta D^2 -\eta D^2 d_{2i-2}^-) \endalign $$
and, setting $d^*=\max\{ d_{2i-2}^-,-d_{2i+5}^+\}$,
$$\min\{s_1^* (2i),s_3^* (2i+3)\} \leq (1-3\beta +D+Dd^*)(1-\eta D+2D^2-\eta D^2-\eta D^2 d^*) \leq \delta_{\infty,4} , $$
since $d^*\leq D/(1-D) -2\beta D/(1-D^2)$. Hence $M^*(\alpha,\gamma)\leq \delta_{\infty,4}$ for $\pm \gamma$ not in $S_{-1}$ or $S_{-4}$ where
$$\delta_{-4}=\delta_{\infty,4}-\frac{2D^2}{(1-D^2)^2} \left(3-3\beta
-3\eta+4D-2D^2+\beta D^2+\eta D^2-D^3\right) <\delta_{\infty,4}$$
and plainly $\delta_{-1} > \delta_{\infty,4}$.

We  are left to  deal with $r=2$, $m\geq 2$.

Suppose first that $b=ma+2$ with $m\geq 3$. Then from above we can assume that 
the expansion of $\gamma$ has $t_{2i+1}=1,$ $t_{2i+2}=-(m-2),-m$ or $-(m+2)$. 
Now if $t_{2i+2}=-(m-2)$ and $t_{2i} $ (or $t_{2i+4}$)
$\geq -m$ then $s_3^*(2i+1)$ is bounded by
$$\align
 & \left(1-2\eta +mD-\frac{\eta D}{1-D}+\frac{(m+2) D^2}{1-D}\right)\left( 1-\beta +(m-2)\beta -\frac{D}{1-D} +\frac{(m+2)\beta D}{1-D} \right)  \\
 & = \delta_{\infty, 6} -\frac{2Dv}{(1-D^2)}(2-3D-3D^2)+\frac{2D^2(1-3D)}{1-D^2}-\frac{2\beta D^2 (1-5D^2)}{(1-D)(1-D^2)} <\delta_{\infty,6}. \endalign $$
Hence we assume that $t_{2i+2}=-(m-2)$ implies that 
$t_{2i+4},t_{2i}=-(m+2)$. So if $\pm \gamma$ is not in $S_0$ or $S_{-3}$ we can assume that $t_{2i}=-(m+2)$, $t_{2i+2}\leq -m$ infinitely often and
$$ \align
s_1^*(2i) &\leq \left(1-\beta -\frac{(m+2)\beta}{1-D^2}+\frac{D}{1-D} -\frac{(m-2)\beta D}{1-D^2}\right)\left(1+\frac{\eta D}{1-D}-\frac{m D}{1-D}\right) \\
 & =\delta_{\infty,6}. \endalign $$
Checking
$$\delta_{0}=\delta_{-3} +\frac{2}{b}\left(1-3\eta +D+\frac{2\eta }{b}\right) >\delta_{-3}, $$
and
$$ \delta_{-3}=\delta_{\infty,6}+\frac{2D}{1+D}\left(1-3\beta -v +\frac{2\beta D
}{1+D}\right)> \delta_{\infty,6} $$
gives the claimed result.

Suppose that $b=2a+2$. We can assume here that $t_{2i+1}=\delta_i, $ $t_{2i+1\pm 1}=0,-2\delta_i$ or $-4\delta_i$ with $\delta_i=\pm 1$. We note the lower bound
$$\delta_{\infty,7}\geq 1-5\beta +D-\eta D. $$
If the expansion contains blocks of the form
$t_{2i}=0$, $t_{2i+1}=-1,$ $t_{2i+2}\neq 4$ then
$$\align
s_1^*(2i)& \leq \left(1-\beta +\frac{D}{1-D}\right)\left(1-\eta -\eta +2D-\frac{\eta D}{1-D} +\frac{4D^2}{1-D}\right) \\
  & < 1-2\eta -\beta +5D-2\beta D-3\eta D+8D^2 \\
  & = 1-5\beta +D-\eta D -6\beta D+8D^2 <\delta_{\infty,7}. \endalign $$
Hence we assume $t_{2i}=0 $ implies that 
$t_{2i-2},t_{2i+2}=\pm 4$. Similarly if $t_{2i+2}=-4$, 
$t_{2i+1}=1$, and
$t_{2i}\leq -2$ then 
$$\align 
s_1^* (2i+1) & \leq \left(1-\eta +\eta -2D +\frac{\eta D}{1-D}\right)\left(1-\beta -4\beta +\frac{D}{1-D}\right) \\
  & < 1-5\beta -D+\eta D +10\beta D-4D^2 \\
 & = 1-5\beta +D-\eta D -2\beta D\left(2a-\frac{2}{a}+\eta -5\right)
 <\delta_{\infty,7}. \endalign $$
So we may assume that $t_{2i}=4\delta_i$ implies that
$t_{2i-2},t_{2i+2}=0$. Hence we need only consider
$\pm \gamma$ in $S_{0}$ or $\gamma$ where the $t_{2i}$ are alternately
$0$ and $-4\delta_i$.
Of these $\gamma$ either $\pm \gamma$ is in $S_{-3}$ or $\pm \gamma$ has infinitely many blocks
$t_{2i},...,t_{2i+4}=-4,1,0,-1,4,$
giving
$$ \align
s_1^* (2i) & = \left(1-\beta +d_{2i}^-\right)\left(1-\eta +\eta -\eta D+\eta Dd_{2i+3}^+\right) \\
s_3^* (2i+3) & = \left(1-\beta-d_{2i+3}^+\right)\left(1-\eta +\eta -\eta D-\eta Dd_{2i}^-\right) \endalign $$
and, setting 
$d^*=\max \{d_{2i}^-,-d_{2i+3}^+\} \leq -\frac{4\beta}{1-D^2} +\frac{D}{1-D} $ 
we obtain
$$ \min\{s^*_1 (2i),s^*_3(2i+3)\}\leq (1-\beta +d^*)(1-\eta D-\eta D d^*)  \leq \delta_{\infty,7}. $$
Finally, $\delta_{0}>\delta_{-3}$ just as for $m\geq 2$ and $\delta_{-3} >\delta_{\infty,7}$ since
$$ \delta_{-3} -\delta_{\infty,7} =\frac{2\beta D}{(1-D)}\left(2-7\eta
 +\frac{2D}{(1+D)(1-D^2)}(11+2D-D^2-3\eta D(1+D))\right). $$

It remains only to deal with  $(a,b)=(3,4),(3,5),(3,6),(5,7),(7,9)$. 

For $(5,7)$ and $(7,9)$ we know from above that $t_{2i+1}=\delta_i$, $\delta_i=\pm 1$, $t_{2i+1\pm 1}=\pm \delta_i$ or 
$-3\delta_{i}$.
Now if $t_{2i}=-3\delta_i$ then we can assume that $t_{2i+2}=t_{2i-2}=\delta_i$, since if for example
$t_{2i}=-3$, $t_{2i+1}=t_{2i-1}=1,$ $t_{2i+2}\leq -1$ then
$$ s_1^*(2i)\leq \left(1-\beta -3\beta +\frac{D+\beta D}{1-D}\right)\left(1-\eta +\eta-D +\frac{\eta D+D^2}{1-D}\right)<\delta_{\infty,8}. $$

Also, if $t_{2i-1}=t_{2i}=t_{2i+1}=-\delta,$ $\delta=\pm 1$, then $t_{2i-2}=t_{2i+2}=3\delta$, since if $t_{2i+1}=t_{2i}=t_{2i-1}=-1$ 
and $t_{2i-2}\leq 1$ then $s_1^* (2i-1)$ is at most
$$ \left( 1-2\eta +\frac{D+\eta D}{1-D}\right)\left(1-2\beta -D +3\beta D-D^2 -\beta D^2+\frac{D^3+\beta D^3}{1-D}\right)
 < \delta_{\infty,8}. $$

Now if $t_{2i+1 }=-\delta$ we can assume that 
we do not have $t_{2i}=t_{2i+2}=-\delta $ since $t_{2i}=t_{2i+1}=t_{2i+2}=-1$ 
would imply
$$ s_1^*(2i)\leq \left(1-\beta -\beta +\frac{D+\beta D}{1-D}\right)\left(1-\eta -\eta  -D +\frac{\eta D+D^2}{1-D}\right)<\delta_{\infty,8}. $$
  
Similarly if $t_{2i+1}=\delta$, $t_{2i+1\pm 1}=-\delta$ and 
$t_{2i-1}$ (or $t_{2i+3}$) $=-\delta$
then $t_{2i-2}$ (or $t_{2i+4}$) $=3\delta$, since for example
$t_{2i-1},t_{2i},t_{2i+1},t_{2i+2}=-1,-1,1,-1$ 
 and $t_{2i-2}\leq 1$ then $s_1^*(2i-1)$ is at most
$$ \left(1-2\eta +D+\eta D-D^2+\frac{\eta D^2+D^3}{1-D}\right)\left( 1-2\beta +D-\beta D+\frac{D^2+\beta D^2}{1-D}\right) < \delta_{\infty,8}. $$

With these exclusions we readily find that $M^*(\alpha,\gamma)<\delta_{\infty,8}$ if the $t_{2i+1}$ are eventually constant and $\pm \gamma$ is not in
$S_0$ or $S_{-3}$ and if the $t_{2i+1}$ eventually alternate and $\pm \gamma$
is not in $S_{-4}$. Hence we may assume that $\pm \gamma$ has infinitely many
blocks of the form $t_{2i-1}=t_{2i+1}=1$ and 
$t_{2i+3}=-1$. Now $t_{2i}=-1$ or $-3$.  If $t_{2i}=-3$ then $s^*_1(2i)$ is at most
$$ \left(1-4\beta +
\frac{D+\beta D +D^2-3\beta D^2}{1-D^2}\right)\left(1 +D-\eta D+\frac{D^2+\eta D^2-3D^3+\eta D^3}{1-D^2}\right) $$ 
which equals $\delta_{\infty,8}.$
So we assume that $t_{2i}=-1$. If $t_{2i+2}=-1$ then 
$t_{2i+4}=3$ and
$$ s^*_3(2i+3)\leq \left( 1+D -\eta D+D^2\right)\left(1-4\beta +\frac{D+\beta D+D^2-3\beta D^2}{1-D^2}\right)<\delta_{\infty,8}. $$
If $t_{2i+2}=1$ then we can assume that $t_{2i+4}=-1$. But if $M^*(\alpha,\gamma)>\delta_{\infty,8}$ then the previous
exclusions force the subsequent coefficients $t_{2i+1+2j}$ to alternate between
$\pm 1$ and the $t_{2i+2j+2}$ between  
$\pm 1$ so that $\gamma $ is in $S_{-4}$. 
One checks numerically that $\delta_{0}>\delta_{-3}>\delta_{-4}>\delta_{\infty,8}.$

Finally we deal with $a=3$, $b=4,5,6$.

In all cases we can rule out $t_{2i}=b$ since $\beta <\delta_{\infty,i}$,
$i=10,11$ or 12.

We begin with $(a,b)=(3,4)$. Suppose that $t_{2i+1}=a$ then 
$$\align  s_2^*(2i) & =\eta (1+\beta +d_{2i}^-)(1-\beta -d_{2i+1}^+),\\
  s_4^*(2i+1) & =\eta (1+\beta +d_{2i+1}^+)(1-\beta -d_{2i}^-), \endalign $$
and if $d_{2i+1}^+\geq aD-\beta D(2+\eta)/(1-D)$ the minimum of these
is certainly at most $\eta (1-(\beta +d_{2i+1}^+)^2)<\delta_{\infty,10}$.
Hence we can assume that $t_{2i+2},t_{2i+3},...$ (and likewise $t_{2i},t_{2i-1},...$) must
take the form $-2,...$ or $0,\pm 1,...$. Also we can successively rule out blocks $t_{2i+1},t_{2i+2},$ of the 
form $-1,-2$ (and hence also
blocks $a,0,1$, since its negative is $a,-2,-1$) or 
$-1,0,-1,0$
or $1,-2,1$ (or
their `negatives') else
$$ s_1^*(2i+1)\leq \left(1-2\eta +\frac{2D+\eta D}{1-D}\right)\left(1-3\beta +\frac{aD}{1-D}\right)<\delta_{\infty,10}, $$
or
$$ s_1^*(2i+2)\leq \left(1-2\eta + aD\eta \right)\left(1-\beta -D +\frac{2\beta D-D^2}{1-D}\right)<\delta_{\infty,10}, $$
or
$$ s_1^*(2i+2)\leq \left(1-3\beta +D+D^2\right)(1+\eta D)<\delta_{\infty,10} . $$
In particular $t_{2i+1}=-1$ forces $t_{2i}=t_{2i+2}=0$. These also
 rule out blocks $-1,0,-1$  and we can further dismiss
$t_{2i-1},t_{2i},t_{2i+1}=1,0,-1$ 
since 
$$\align
s_1^*(2i+1) & = \eta (1-\beta +D - D(-d_{2i-2}^-))(1-\beta +d_{2i+1}^+) \\
s_3^*(2i-2) & =  \eta (1-\beta +D - Dd_{2i+1}^+)(1-\beta +(-d_{2i-2}^-)) 
\endalign $$
and since $d_{2i+2}^+\leq (aD-D^2)/(1-D^2)$ the minimum of these is
certainly at most $\eta (1-\beta+d_{2k+1}^+)(1-\beta +D-Dd_{2k+1}^+)\leq \delta_{\infty,10}$. In particular this rules out blocks $1,0$. Hence we can assume that the sequence for $\gamma$ consists entirely
of blocks $F_2=(a,-2)$, $M=(a,-2,1,-2)$ or $M'=(a,0,-1,0)$. We suppose that $\gamma$ does not have period
$F_2$ or $M$ or $M'$. If a block $t_{2i+1},t_{2i+2},...=M,...$ is preceeded by blocks
$M'$ or $M'F_2$ or $F_2F_2$ then
$$\align d_{2i+1}^+ & \leq -2\beta+D-2\beta D +\frac{aD^2-D^3}{1-D^2},\\
 d_{2i}^- & \geq -2\beta +aD-2\beta D+\frac{aD^2-2\beta D^2+D^3-2\beta D^3}{1-D^2}\endalign $$
and
$s_4^*(2i+1)=\eta (1+\beta+d_{2i+1}^+)(1-\beta-d_{2i}^-)<\delta_{\infty,10}$.
Hence we can assume that the sequence consists of blocks of $M$'s
separated by single $F_2$'s (or its negative with $M'$'s replacing the $M$'s).
Suppose that we have a block $...|t_{2i+1},t_{2i+2},...$ of the form
$...F_2 M^s F_2 | M^t F_2,... $ with $t\geq s+1$ then 
$$ d_{2i+1}^+ \leq \frac{-2\beta +D-2\beta D +aD^2}{1-D^2} +(a-1)D^{1+2t}, $$
$$ d_{2i}^-\geq -2\beta +\frac{aD-2\beta D+D^2-2\beta D^2}{1-D^2}+(a-1)D^{2+2s}(1-D)$$
and $s_4^*(2i+1)\leq \delta_{\infty,10}+\eta (a-1)D^{2s+2}(D(1+\beta D)-(1-\beta)(1-D))<\delta_{\infty,10}$. This leaves only the periodic elements of $S_{k,10}$.

Suppose next that $(a,b)=(3,5)$. If $t_{2i+1}=a$
then in the same way we can rule out $d_{2i+1}^+\geq (\beta-D)/(1-D)$
or $d_{2i+1}^+\leq -3\beta +(D-\beta D)$ and deduce that $t_{2i+2},...$ (and $t_{2i},...$)
take the form $1,-1,...$ or $-1,...$ or $-3,a,...$ or $-3,1,...$. We can also successively rule out blocks 
$t_{2i},t_{2i+1},$ of the form $-3,-1$
and $t_{2i-1},t_{2i},...=1,-3,1,...$ else
$$ s_1^*(2i)\leq \left(1-4\beta +\frac{aD+\beta D}{1-D}\right)\left(1-2\eta +\frac{3D+\eta D}{1-D}\right)<\delta_{\infty,11} $$
or
$$ s_1^* (2i) \leq \left(1-4\beta +\frac{D+\beta D}{1-D}\right)\left( 1+\frac{D+\eta D}{1-D}\right)<\delta_{\infty,11}. $$
In particular there are no $t_{2i}=3$ and the $t_{2i}=-3$ are confined to blocks $a,-3,a$ or
$a,-3,1$ or $1,-3,a$.
Likewise we can rule out blocks
$t_{2i},t_{2i+1}=-1,-1$ with $t_{2i+2}=-1$ or with $t_{2i-1}\leq 1$ else
$$ s_{1}^*(2i)\leq (1-2\beta +aD +\beta D)(1-2\eta -D+a\eta D+D^2)<\delta_{\infty,11} $$
or
$$ s_{1}^*(2i)\leq (1-2\beta +D +\beta D+D^2)(1-2\eta +D+3\eta D+D^2)<
\delta_{\infty,11} $$
respectively (in particular this rules out blocks
$1,1$). Hence we can successively
rule out $t_{2i},t_{2i+1}=-3,a$ with $t_{2i+2}=1$ or with $t_{2i+2}=-1$ else
$$ s_2^*(2i)\leq \eta (1+\beta -3\beta +aD+\beta D)(1-\beta -\beta +D+\beta D)
<\delta_{\infty,11}, $$
$$ s_2^*(2i)\leq \eta (1+\beta -3\beta +aD)(1-\beta +\beta +D)<\delta_{\infty,11}, $$
and blocks $t_{2i+1},...=a,-3,a,-3,1$
else
$$ s_4^*(2i+3)\leq \eta (1+\beta -3\beta +D)(1-(\beta-3\beta +aD-3\beta D))<\delta_{\infty,11}. $$
Hence either $\gamma$ has period $F_1=(a,-1)$ or period $F_3=(a,-3)$ or the 
$t_{2i-1}=a$ 
occur inside blocks $H=(1,-3,a,-3,1)$, $H'=(-1,1,a,1,-1)$ or $\pm 1,-1,a,-1,\pm 1$.
Moreover any block $H$ or $H'$ must be followed (and preceded) by a block $G=(-1,a,-1)$ since
$t_{2i-1},...,t_{2i+3}=a,-3,1,-1,1$  
leads to $\min\{s_2^*(2i),s_4^*(2i+1)\}\leq \eta (1-(\beta -3\beta+D-\beta D+\frac{D^2-\beta D^2+a D^3-\beta D^3}{1-D^2})^2)$ and any $G$ must be followed  by an $H$ or $H'$ since blocks
$t_{2i},...,t_{2i+3}=-1,-1,1,-1$ 
give
$$ s_1^*(2i)\leq (1-2\beta +aD -\beta D+D^2)(1-2\eta +D -\eta D+D^2)<\delta_{\infty,11}. $$
Hence apart from period $B_1$ (where $\delta_0<\delta_{\infty,11}$) 
or $F_2$ (or their negatives)
we can assume the sequence consists of blocks $HG$ or $H'G$.
We suppose that $\gamma$ is not in $S_{-9}$ so that we have
blocks $...|t_{2i-1}...=...|HGHG...$ and write 
$$ d_{2i}^-=-3\beta +D -\lambda +D^3w^- ,\;\;\;\; d_{2i+1}^+= -3\beta +D+  \lambda+ D^7 w^+,$$
where $\lambda :=\frac{D^3(1-2\beta-2\beta D+D^2)}{1+D^4}$, so that 
$$\align s_{2}^*(2i) & =\eta \left(1-2\beta +D - \lambda + D^3w^-\right)(1+2\beta -D -\lambda - D^7 w^+) \\
s_{4}^*(2i+9) & =\eta \left(1-2\beta +D - \lambda + D^3w^+\right)(1+2\beta -D -\lambda - D^7 w^-). \endalign $$
Hence if $w^+$ and $w^-$ are of opposite signs then one of these is certainly
less than $\delta_{\infty,11}$ and if both are negative (say $w^-<w^+$)
then $s_2^*(2i)\leq \delta_{\infty,11}+D^3w^-(1-D^4)<\delta_{\infty,11}$.
Hence we can assume that $w^-,w^+>0$ and hence that we are in a block
$...HG (HGH'G)^l HGHG(H'GHG)^k HG...$. Now 
$$w^+\geq 2(1-2\beta-2\beta D+D^2)\frac{D^{8k}(1-D^4-D^8)}{(1-D^8)},\;\;w^-\leq 2(1-2\beta-2\beta D+D^2)\frac{D^{8l}}{(1-D^8)} $$
hence if $l\geq k+1$ we have 
$$ s_2^*(2i)\leq \delta_{\infty,11}-2(1-2\beta-2\beta D+D^2)\frac{D^{8k+7}}{1-D^8}\left((1-D^4-D^8)(1-2\beta)-D^4 (1+2\beta)\right)$$
less than $\delta_{\infty,11}$, leaving only the elements of $\delta_{k,11}$.


Suppose finally that $(a,b)=(3,6)$. Now if $t_{2i+1}=a$ then (in the same way as before) we can assume that $-2\beta <d_{2i}^-,d_{2i+1}^+<0$ and hence that
$t_{2i},...$ and $t_{2i+2},...$ are of the form $0,-1,...$ or $-2,1,...$ or  $-2,a,...$.
Now we can rule out $t_{2i}=-4$ else
$$s_1^*(2i)\leq \left(1-5\beta +\frac{D+4\beta D}{1-D}\right)\left(1-\eta +\frac{\eta +4D}{1-D}\right)<\delta_{\infty,12}$$
and blocks $t_{2i+1},t_{2i}=-1,-2$ else
$$ s_1(2i)\leq \left(1-2\eta +\frac{2D+\eta D}{1-D}\right)\left(1-3\beta +\frac{D+2\beta D}{1-D}\right)<\delta_{\infty,12} $$
and $t_{2i+1},t_{2i+2},t_{2i+3}=-1,0,-1,$ else
$$ s_1^*(2i+1)\leq \left(1-2\eta +\frac{2D-\eta D}{1-D}\right)\left(1-\beta -\frac{(D-2\beta D)}{1-D}\right)<\delta_{\infty,12}. $$
We can also rule out blocks $t_{2i+1},t_{2i+2},t_{2i+3}=-1,0,1$ else
$$\align s_1^*(2i+1)& =(1-2\eta +\eta d_{2i}^-)(1-\beta +D-D(-d_{2i+2}^+)) \\
 s_3^*(2i+2)& =(1-2\eta +\eta (-d_{2i+2}^+))(1-\beta +D-Dd_{2i}^-)) \endalign $$
with $d_{2i}^-\leq (2\beta -D)/(1-D)$ and the smallest of
these is at most $(1-2\eta+\eta d_{2i}^-)(1-\beta +D-Dd_{2i}^-)\leq \delta_{\infty,12}$. In particular this means that $t_{2i+1}=1$
occur only in blocks $-2,1,-2$. Now if $t_{2i+1},t_{2i+2},t_{2i+3}=a,-2,1$ then $t_{2i},t_{2i-1}=-2,1$ else
$d_{2i+1}^+\leq -2\beta +D$ and $d_{2i}^-\geq -2\beta +aD-2\beta D$ and $s_4^*(2i)\leq \eta (1+\beta D)(1-\beta +D)<\delta_{\infty,12}$. Hence either $\gamma$ or its negative  has period $B_2'=(1,-2)$ or $F_2=(a,-2)$
or takes the form
$...|t_{2i+1}...=...F_2 {B_2'}^l | F_2 {B_2'}^k F_2...$
with $d_{2i}^-,d_{2i+1}^+=(-2\beta+D)/(1-D)+w^-,w^+$ where $w^+\geq (a-1)D^{k+1}$, $w^-\leq (a-1)D^{l+1}/(1-D)$ and if $l\geq k+1$ then
$s_2^*(2i)\leq \delta_{\infty,12}-\eta (a-1)D^{k+1} ((1-\beta)-D(1+\beta)/(1-D))$. This leaves only the elements of $S_{k,12}$. $\blacksquare$

\head 6.  Proof of Theorem 3 ($a=2$)  \endhead      

We assume that $\alpha$ has $a_{2i-1}=a=2$, $a_{2i}=b\geq 5$, and that $\gamma$ has $M^*(\alpha,\gamma) >\delta_{\infty}$.
If $t_{2i-1}=0$ then
$$\align 
s_1^*(2i-1) & =\eta (1-\beta+d_{2i-2}^-)(1-\beta +d_{2i-1}^+) \\
s_3^*(2i-1) & =\eta (1-\beta-d_{2i-2}^-)(1-\beta -d_{2i-1}^+) \endalign $$
and if $d=\max\{\left|d_{2i-2}^-\right|,\left|d_{2i-1}^+\right|\}$
then the minimum of these is certainly at most $\eta ((1-\beta)^2-d^2)$.
Hence when $b$ is even we must have $d\leq \beta$ ruling out $t_{2i-1\pm 1}\geq 2$
or $t_{2i-1\pm 1}\leq -4$ or $t_{2i-1\pm 1}=-2,t_{2i-1\pm 2}=0$ and hence $t_{2i-1\pm 1},t_{2i-1\pm 2},...=0,...$ or $-2,a,...$. Likewise if $b$ is odd we must have
$d\leq \frac{1}{2}\beta^2$ ruling out $t_{2i-1\pm 1}\geq 1$ or $t_{2i-1\pm 1}\leq -3$ or
$t_{2i-1\pm 1},t_{2i-1\pm 2}=-1,0$ and any $t_{2i-1}=0$ must occur inside a block
$t_{2i-3},...,t_{2i+1}=a,-1,0,-1,a$.

If $t_{2i-1}=a$ then
$$\align 
s_2^*(2i-2) & =\eta (1+\beta+d_{2i-2}^-)(1-\beta -d_{2i-1}^+) \\
s_4^*(2i-1) & =\eta (1-\beta-d_{2i-2}^-)(1+\beta +d_{2i-1}^+). \endalign $$
We suppose that $t_{2j}\leq t$ for all $j$ with $t_{2i}=t$ infinitely often.
If $t\geq 1$ then in view of the above we can assume that they occur inside blocks
$t_{2i-1},t_{2i},t_{2i+1}=a,t,a$ and hence that $t\leq b-4$.
Since
$$\align 
s_2^*(2i-2) & =\eta (1-\beta+(2\beta+d_{2i-2}^-))(1-\beta -(t\beta +aD)+2\beta D -D(2\beta+d_{2i+1}^+)) \\
s_4^*(2i+1) & =\eta   (1-\beta+(2\beta+d_{2i+1}^+))(1-\beta -(t\beta +aD)+2\beta D -D(2\beta+d_{2i-2}^-))     . \endalign $$
and if $d=\min \{ 2\beta +d_{2i-2}^-,2\beta +d_{2i+1}^+\}$ then certainly
$$ \min\{s_2^*(2i-2),s_4^*(2i+1)\}\leq  \eta (1-\beta+d)(1-\beta -(t\beta +aD)+2\beta D -Dd). $$  
Now if $t_{2i-2}$ or $t_{2i+2}\leq t-2$ then $d\leq (t\beta +aD)/(1-D)$ and, as $(t\beta +a D)/(1-D)\leq 1 \leq (b-t-3+D)/(1+D)=\frac{1}{2}\left((1-\beta-(t\beta+aD)+2\beta D)/D -(1-\beta)\right)$, 
$$\align \min\{s_2^*(2i-2),s_4^*(2i+1)\}& 
\leq  \eta \left( (1-\beta +\beta D)^2 -\left(\frac{t\beta +aD}{1-D}-\beta D\right)^2 \right) \\
 &\leq  \eta ((1-\beta)^2+2\beta D-((t+1)\beta)^2) <\delta_{\infty}.\endalign  $$
Hence either $\gamma$ has period $(a,t)$ or  $t\leq 0$. We assume that
$\gamma$ is not in $S_{0,t}$ (these were dealt with in Lemma 3). 
Now observing that if $\gamma$ has blocks $(t_{2i-1},t_{2i},t_{2i+1})=(0,-\lambda,0)$ or $(a,-\lambda,0)$
or $(a,-\lambda,a)$ then $1-\alpha-\gamma$  
contains corresponding blocks $(0,\lambda,0)$ or $(a,\lambda -2,0)$
or $(a,\lambda -4,a)$  we readily deduce that our 
sequence consists of blocks $(a,-1,0,-1,a)$ or $(a,-1,a)$ or $(a,-3,a)$ 
when $b$ is odd, and when $b$ is even that $t_{2i-1}=a$ are adjoined by
0 or $-2$ or $(-4,a)$. Now if $b$ is even and  $t_{2i-3},...,t_{2i-1}=a,-2,0$
and $t_{2i}=-2$ or $t_{2i},t_{2i+1}=0,0$ then
$$ \align s_2^*(2i-1) & \leq \eta \left(1-\beta -2\beta +\frac{aD}{1-D}\right)\left(1-\beta +\frac{aD^2}{1-D}\right)\\
 &  \leq \eta \left( (1-\beta)^2-\beta (1-\beta -3D)\right)<\delta_{\infty}.
\endalign $$
 Hence when $b$ is even either $\gamma$ has period $0,0$ or the $t_{2i-1}=0$
occur inside blocks $a,0,0,-2,a$ or their negatives $a,-2,0,0,a$.
Hence we can assume that the expansion of $\gamma$ consists
entirely of blocks $(a,0,0,-2)$ or $(a,-2,0,0)$ or $(a,-2)$ or $(a,0)$ or $(a,-4)$ when $b$ is even
and blocks $(a,-1,0,-1)$ or $(a,-1)$ or $(a,-3)$ when $b$ is odd.

We deal first with $b$ even. Suppose that we have a block $t_{2i-1},t_{2i},t_{2i+1}=a,0,a$ then writing $d_{2i-2}^-,d_{2i+1}^+=-\beta +w^-,-\beta+w^+$ (where $-\beta =(-2\beta+aD)/(1-D)$) we have
$$\align s_2^*(2i-2) & =\eta (1+w^-)(1-2\beta-Dw^+) \\
s_4^*(2i+1) & =\eta (1+w^+)(1-2\beta-Dw^-) \endalign $$
Now writing $w=\min\{w^+,w^-\}$ the minimum of these is plainly
at most $\eta (1+w)(1-2\beta -Dw)$ and we can assume that $w^+,w^->0$.
So if $t_{2i-1},...,t_{2i+1}=a,0,a$ then the adjacent
$t_{2i-2},t_{2i-3},...$ or  $t_{2i+2},t_{2i+3},...$
take the form 
$(-2,a)^l\{(0,a) \hbox{ or } (0,0,-2,a)\},...$,
and if  $t_{2i-1},...,t_{2i+1}=a,-4,a$, 
take the form $(-2,a)^l\{(-4,a) \hbox{ or } (-2,0,0,a)\},...$ .

Similarly suppose that we have a block $t_{2i-3},...,t_{2i+1}=a,0,0,-2,a$ 
then writing $d_{2i-4}^-,d_{2i+1}^+=-\beta +v^-,-\beta +v^+$
then 
$$\align s_3^*(2i-2) & =\eta (1-Dv^+)(1-2\beta-Dv^-) \\
s_1^*(2i-1) & =\eta (1+Dv^-)(1-2\beta+Dv^+) \endalign $$
the previous conditions $|d_{2i-2}^-|,|d_{2i-1}^+|\leq \beta$ 
giving $v^+\geq 0$ and $v^-\leq 0$. Hence this block must be 
contained inside a block
$$ \{(a,-4) \hbox{ or } (a,0,0,-2) (a,-2)^m (a,0,0,-2,a) (-2,a)^l \{(0,0,-2,a) \hbox{ or } (0,a)\} $$ 
Combining these  two restrictions it is clear that $\gamma$ 
(it or its negative) must  consist solely of blocks of the form $(a,0)(a,-2)^{l_i}$
or solely of the form $(a,0,0,-2)(a,-2)^{l_i}$. 

Suppose then that we have a block
$$ ...(a,0)(a,-2)^m (a,0,a) (-2,a)^l,(0,a),... $$
with $m\geq l+1$. Then $w^+\geq 2\beta D^l $ and $w^-\leq 2\beta D^m/(1-D)$
and
$$\align  s_2^*(2i-2) & \leq \eta \left(1+\frac{2\beta D^{l+1}}{1-D}\right)(1-2\beta -2\beta D^{l+1}) \\
   & \leq \eta \left( 1-2\beta -\frac{2\beta D^{l+1}}{1-D}((1-D)-(1-2\beta))\right)<\delta_{\infty}. \endalign $$
Likewise if we have blocks
$$ ...(a,0,0,-2) (a,-2)^m (a,0,0,-2,a)(-2,a)^l (0,0,-2,a)... $$
with $l\geq m+1$ then
$$ v^-\leq -(aD-2\beta D)D^m=-\beta (1-D)D^{m},\;\;\; v^+\leq (2\beta -aD)D^l/(1-D)\leq \beta D^{m+1} $$
and
$$ s_1^*(2i-1)\leq \eta \left( 1-2\beta -\beta D^{m+1}((1-D)(1-2\beta)-D)\right)<\delta_{\infty}. $$
This leaves only the periodic elements of $S_{2k}$ and $S_{2k+1}$.

Suppose now that $b$ is odd and that we have a block
$t_{2i-3},...,t_{2i+1}=a,-1,0,-1,a$. We assume without loss of generality
that $d_{2i+1}^+\leq d_{2i-4}^-$ (else reverse their roles) and
write
$$ d_{2i-4}^-=\frac{-\beta +aD-3\beta D+aD^2}{1-D^2}+w^-,\;\;\;d_{2i+1}^+=\frac{-3\beta +aD -\beta D+aD^2}{1-D^2}+w^+, $$ 
so that
$$\align 
s_1^*(2i-1) & = \eta \left(1-\beta -\frac{1}{2}\beta^2 +Dw^+\right)\left(1-\beta +\frac{1}{2}\beta^2+Dw^-\right), \\
s_3^*(2i-1) & = \eta \left(1-\beta -\frac{1}{2}\beta^2 -Dw^-\right)\left(1-\beta +\frac{1}{2}\beta^2-Dw^+\right), \endalign $$
the previous restriction $|d_{2i-2}^-|,|d_{2i-1}^+|\leq \frac{1}{2}\beta^2$
forcing $-2\beta +\frac{1}{2}\beta^2\leq d_{2i-4}^-,d_{2i+1}^+\leq -\frac{1}{2}\beta^2$ and hence $w^-\leq 0$, $w^+ \geq 0$. So the sequence must
take the form 
$$\{ (a,-3) \hbox{ or } (a,-1,0,-1) \} (a,-3,a,-1)^l (a,-1,0,-1,a) (-3,a,-1,a)^m \{ (-1,a) \hbox{ or } (-1,0,-1,a)\} $$
(or its reflection) with
$$ w^-\leq \left(-aD +\frac{2\beta D}{1-D^2}\right) D^{2l},\;\;\;\; w^+\leq \frac{2\beta D^{2m}}{1-D^2}. $$
Hence if $m\geq l+1$ we have
$$ s_1^*(2i-1)\leq \delta_{\infty}+ \frac{ D^{2l+2}}{1-D^2} \left(2\left(1-\beta+\frac{1}{2}\beta^2\right)D^2-(1-D-D^2-D^3)\left(1-\beta-\frac{1}{2}\beta^2\right)\right)<\delta_{\infty}.  $$
Taking negatives similarly dismisses $l\geq m+1$ so we can assume that $l=m$ for such blocks. 

Suppose now that we have a block $t_{2i-3},...,t_{2i+1}=a,-1,a,-1,a$
then writing
$$ d_{2i-4}^-,d_{2i+1}^+=\frac{-3\beta +aD-\beta D+aD^2}{1-D^2} +v^-,v^+$$
we have
$$\align 
s_2^*(2i-4) & = \eta \left(1-\beta+\frac{1}{2}\beta^2+v^-\right)\left(1-\beta -\frac{1}{2}\beta^2-D^2v^+\right) \\
s_4^*(2i+1) & = \eta \left(1-\beta+\frac{1}{2}\beta^2+v^+\right)\left(1-\beta -\frac{1}{2}\beta^2-D^2v^-\right). \endalign $$
Now we can assume that $v^-,v^+>0$ else $v=\min\{v^+,v^-\}<0$ and
the minimum of these is at most $\eta (1-\beta+\frac{1}{2}\beta^2+v)(1-\beta-\frac{1}{2}\beta^2-D^2v)<\delta_{\infty}$.
In particular we can rule out blocks $(a,-1,0,-1,a)(-3,a,-1,a)^k(-1,a)$
or $(a,-3)(a,-3,a,-1)^k (a,-1,0,-1,a)$  when $k\geq 1$. Hence we can assume that the blocks $(a,-1,0,-1,a)$ only occur in the 
periodic elements of $S_{2k+1}$. Since we have also ruled out blocks
$(a,-3)$ $(a,-3,a,-1)^k$ $(a,-1)$, $k\geq 1$, it remains only to consider  $\gamma$ 
in $S_{-2}$ or $\gamma$ (or its negative) consisting
solely  of blocks $(-1,a)(-3,a,-1,a)^{l_i}$ with the $l_i\geq 0$.
Suppose then that we have blocks
$$ ....(a,-1)(a,-1,a,-3)^l (a,-1,a,-1,a) (-3,a,-1,a)^m (-1,a) (-3,a,-1,a)^k (-1,a)...  .$$
Now 
$$v^- \geq \left(2\beta -\frac{2\beta D}{(1+D)}\right)D^{2l} =\frac{2\beta D^{2l}(1-D)}{1-D^2}, $$ 
and  if $m\geq l+2$
$$ v^+\leq \frac{2\beta D^{2m}}{1-D^2} \leq \frac{2\beta D^{2l+4}}{1-D^2}, $$
and if $m=l+1$, $k\geq 1$,
$$ v^+\leq \left(2\beta -2\beta D+\frac{2\beta D^2}{1-D^2} \right)D^{2m}, $$
and in either case 
$$ v^+\leq \frac{2\beta D^{2l+2}(1-D+D^3)}{1-D^2},$$
giving
$$ s_4^*(2i+1)\leq \delta_{\infty} +\frac{2 D^{2l+3}}{1-D^2}\left(
(1-D+D^3)\left(1-\beta-\frac{1}{2}\beta^2\right)-(1-D)\left(1-\beta+\frac{1}{2}\beta^2\right)\right)< \delta_{\infty}. $$ 
Hence  $m\leq l+1$ and if $m=l+1$ we must have $k=0$
and (since similar reasoning gives $m\leq k+1$) either $l=m$ or  $(l,m,k)=(0,1,0)$. 
But in this latter case looking at the preceding block
$$...(a,-1) (a,-1,a,-3)^{r}(a,-1)(a,-1,a,-1,a)(-3,a,-1,a)(-1,a)(-1,a),....$$
we must have $r=0$ or 1. But if $r=0$ then $v^-\geq 2\beta$, $v^+\leq 2\beta D^2/(1-D^2)$ and 
$$ s_4^*(2i+i)\leq \delta_{\infty} + \frac{2 D^3}{1-D^2}\left((1-\beta-\frac{1}{2}\beta^2)-(1-D^2)(1-\beta+\frac{1}{2}\beta^2)\right)<\delta_{\infty}. $$
Hence this case only occurs when we alternate between blocks $(-1,a)(-3,a,-1,a)^l$ of length $l=1$ and
$l=0$ giving the elements of $S_{-1}$. The remaining $\gamma$ have blocks
of equal length giving the elements of $S_{2k}$. $\blacksquare$


\Refs
 
\ref{Barnes} \key{\bf 1} \by   E. S. Barnes  \& H. P. F. Swinnerton-Dyer,
{\it  The inhomogeneous minima of binary quadratic forms.}
Part I, {\it Acta Math. } {\bf 87} (1952), 259-323, Part II, {\it Acta M
ath. } {\bf 88} (1952), 279-316, Part III, {\it Acta Math. } {\bf 92} (1
954), 199-234,
Part IV (without second author) {\it Acta Math. } {\bf 92} (1954), 235-2
64.
 \endref 
 
\ref \key{\bf 2} \by  D.\ Cardon  \paper  
A Euclidean ring containing $\Bbb Z [\sqrt{14}]$
 \jour C. R. Math. Rep. Acad. Sci. Canada   \vol 19  \yr 1997  \pages  28--32 \endref



\ref \key{\bf 3} \by  D. J. Crisp, W. Moran \& A. D. Pollington  \paper  
Inhomogeneous diophantine approximation and Hall's ray
 \it preprint  \endref
   
\ref \key{\bf 4} \by  T. W. Cusick and M. E. Flahive \book  
The Markoff and Lagrange spectra,
 mathematical surveys and monographs, no. 30, American Mathematical Society \yr 1989 \endref



 

 
\ref \key{\bf 5} \by  H. J. Godwin  \paper  
On a conjecture of Barnes and Swinnerton-Dyer
 \jour Proc. Camb. Phil. Soc.  \vol 59 \yr 1963 \pages 519--522 \endref

\ref \key{\bf 2} \by  M.\ Harper  \paper  
Reference??
 \jour   \vol  \yr  \pages  \endref

\ref \key{\bf 5} \by  F. Lemmermeyer   \paper  
The Euclidean algorithm in algebraic number fields
 \jour Expo. Math.  \vol 13 \yr 1995 \pages 385--416 \endref

\ref \key{\bf 6} \by  C. G. Pinner  \paper  
More on inhomogeneous Diophantine approximation
 \jour preprint  \endref

\ref \key{\bf 7}  \by P. Varnavides  \paper  
Non-homogeneous quadratic forms, I 
\jour Proc. K. Ned. Akad. Wet. (Indag. Math.) \vol 51 \yr 1948 
\pages 396--404, 470--481 \endref


 
\endRefs
\enddocument